\input ./style/arxiv-general.cfg
\documentclass[aos,MSNbibl,seceqn,nameyear,dvips]{arximspdf}
\makeatletter
   \@ifpackageloaded{graphicx}{}{\usepackage{graphicx}}
\makeatother
\usepackage{mathbh}
\usepackage[ruled,vlined]{algorithm2e}
\SetKwHangingKw{KwParameters}{parameters}
\SetKwHangingKw{KwInit}{initialization}
\SetKwHangingKw{KwIn}{input}
\SetKwFor{KwAndWhen}{and when}{do}{endfor}

\doi{10.1214/15-AOS1345}
\volume{43}
\issue{6}
\pubyear{2015}
\firstpage{2412}
\lastpage{2450}
\docsubty{FLA}

\makeatletter
\let\my@algocf@latexcaption\algocf@latexcaption
\let\my@addcontentsline\addcontentsline
\long\def\algocf@latexcaption#1[#2]#3{%
\def\addcontentsline##1##2##3{}%
\my@algocf@latexcaption{#1}[#2]{#3}%
\global\let\addcontentsline\my@addcontentsline%
}
\renewcommand{\j}{\mathrm{J}}
\newcommand{\rrvert}{\vert}
\newcommand{\llvert}{\vert}
\renewcommand{\mid}{|}
\newtheorem{theorem}{Theorem}[section]
\newtheorem{proposition}{Proposition}
\newtheorem{corollary}{Corollary}
\newtheorem{lemma}{Lemma}
\newproclaim{definition}{Definition}
\newproclaim{remark}{Remark}
\newproclaim{example}{Example}
\newcommand{\1}{\mathbh{1}}
\newcommand{\rme}{\mathrm{e}}
\newcommand{\N}{\mathbb{N}}
\newcommand{\Z}{\mathbb{Z}}
\newcommand{\R}{\mathbb{R}}
\newcommand{\E}{\mathbb{E}}
\newcommand{\PP}{\mathbb{P}}
\newcommand{\btheta}{\bolds\theta}
\newcommand{\rmd}{\mathrm{d}}
\def\m{m}
\def\mphi{\phi}
\makeatother

\begin{document}
\begin{frontmatter}

\title{Aggregation of predictors for nonstationary sub-linear
processes and online adaptive forecasting of time varying
autoregressive~processes\thanksref{T1}}
\runtitle{Aggregation of predictors}

\begin{aug}
\author[A]{\fnms{Christophe}~\snm{Giraud}\thanksref{M1}\ead[label=e1]{christophe.giraud@math.u-psud.fr}},
\author[B]{\fnms{Fran\c{c}ois}~\snm{Roueff}\thanksref{M2}\ead[label=e2]{francois.roueff@telecom-paristech.fr}}\\
\and
\author[B]{\fnms{Andres}~\snm{Sanchez-Perez}\corref{}\thanksref{M2}\ead[label=e3]{andres.sanchez-perez@telecom-paristech.fr}}
\runauthor{C. Giraud, F. Roueff and A. Sanchez-Perez}
\affiliation{Universit\'e Paris Sud\thanksmark{M1} and T\'el\'ecom
ParisTech; CNRS LTCI\thanksmark{M2}}
\address{C. Giraud\\
D\'epartement de Math\'ematiques\\
B\^atiment 425, Facult\'e des Sciences d'Orsay\\
Universit\'e Paris-Sud\\
F-91405 Orsay Cedex\\
France\\
\printead{e1}}
\address[B]{F. Roueff\\
A. Sanchez-Perez\\
Institut Mines-T\'el\'ecom\\
T\'el\'ecom ParisTech\\
CNRS LTCI\\
37 rue Dareau\\
75014 Paris\\
France\\
\printead{e2}\\
\phantom{E-mail: }\printead*{e3}}
\end{aug}
\thankstext{T1}{Supported in part by the Conseil r\'egional
d'\^{I}le-de-France under a doctoral allowance of its program R\'eseau de
Recherche Doctoral en Math\'ematiques de l'\^{I}le de France (RDM-IdF)
for the period 2012--2015 and by the Labex LMH (ANR-11-IDEX-003-02).}

%
\received{\smonth{5} \syear{2014}}
%
\revised{\smonth{5} \syear{2015}}

%
\begin{abstract}
In this work, we study the problem of aggregating a finite number of
predictors for nonstationary sub-linear processes. We provide oracle
inequalities relying essentially on three ingredients: (1)~a uniform
bound of
the $\ell^1$ norm of the time varying sub-linear coefficients, (2)~a Lipschitz
assumption on the predictors and (3)~moment conditions on the noise appearing
in the linear representation. Two kinds of aggregations are considered giving
rise to different moment conditions on the noise and more or less sharp
oracle inequalities. We apply
this approach for deriving an adaptive predictor for locally stationary time
varying autoregressive (TVAR) processes. It is obtained by aggregating a
finite number of well chosen predictors, each of them enjoying an optimal
minimax convergence rate under specific smoothness conditions on the TVAR
coefficients. We show that the obtained aggregated predictor achieves a
minimax rate while adapting to the unknown smoothness. To prove
this result, a lower bound is established for the minimax rate of the
prediction risk for the TVAR process. Numerical experiments complete this
study. An important feature of this approach is that the aggregated predictor
can be computed recursively and is thus applicable in an online prediction
context.
\end{abstract}

%
\end{frontmatter}

\section{Introduction}\label{sec1}
In many applications where high frequency data are observed, we wish to
predict the next values of this time series through an online prediction
learning algorithm able to process a large amount of data.
The classical stationarity assumption on the distribution of the
observations has to be weakened to take
into account some smooth evolution of the environment. From a
statistical modelling point of view, this is described by some time
varying parameters. In order to sequentially track them from
high-frequency data, the
algorithms must require few operations and a low storage capacity to
update the parameters estimation and the prediction after each new observation.
The most
common online methods are least mean squares (LMS), normalized least
mean squares (NLMS), regularized least squares (RLS) or Kalman. All of
them rely on the choice
of a gradient step, a forgetting factor, or more generally on a tuning
parameter corresponding to some {a priori} knowledge on how
smoothly the local statistical distribution of the data evolves along
the time. To adapt
automatically to this smoothness, usually unknown in practice, we
propose to use an exponentially weighted aggregation of several such
predictors, with various tuning
parameters. We emphasize that to meet the online constraint, we cannot
use methods
that require a large amount of computations (such as cross validation).

The exponential weighting technique in aggregation have been developed
in parallel
in the machine learning community [see the seminal paper \citet{Vovk1990}],
in the statistical community [see \citet{Catoni1997},
\citeauthor{Yang2000} (\citeyear{Yang2000,Yang2004}), \citet{LeungBarron2006}, or more recently
\citet{DalalyanTsybakov2008,Audibert2009,RigolletTsybakov2012}]
and in the game theory community for individual sequences prediction
[see \citet{Cesa-BianchiLugosi2006} and \citet{Stoltz2011} for
recent surveys].
In contrast to the classical statistical setting, in the individual
sequence setting the observations are not assumed to be generated by an
underlying stochastic
process. The link between both settings has been analyzed in \citet
{Gerchinovitz2011} for the regression model with fixed and random designs.

Exponential weighting has also been investigated in the case of weakly dependent
stationary data in \citet{AlquierWintenberger2012}. More recently,
an approach
inspired from individual sequences prediction has been studied in
\citet
{AnavaHazanMannorShamir2013} for
bounded ARMA processes under some specific conditions on the (constant) ARMA
coefficients.

In this contribution, we consider two possible aggregation schemes
based on exponential
weights which can be computed recursively. We provide oracle inequalities
applying to the aggregated predictor under the following main
assumptions that (1) the
observations are sub-linearly with respect to a sequence of random
variables with possibly
time varying linear coefficients and (2) the predictors to be
aggregated are
Lipschitz functions of the past. An important feature of our
observation model
is that it embeds the well-known class of \emph{locally stationary}
processes. We refer to \citet{Dahlhaus2009} and the
references therein for a recent general view about statistical
inference for
locally stationary processes. As an application, we focus on a particular
locally stationary model, that of the time varying autoregressive (TVAR)
process. The minimax rate of certain recursive estimators of the TVAR
coefficients is studied in \citet{MoulinesPriouretRoueff2005}. To our
knowledge, there is not a well-established method on the automatic
choice of
the gradient step when the smoothness index is unknown. Here, we are
interested in the prediction problem which is closely
related to the estimation problem. We show that the
proposed aggregation
methods provide a solution to this question, in the sense that they
give rise to
recursive adaptive minimax predictors.

The paper is organized as follows. In Section~\ref{sec:online-aggr-pred}, we
provide oracle inequalities for the aggregated predictors under general
conditions applying to nonstationary sub-linear processes. TVAR
processes are
introduced in Section~\ref{section:time-varying_autoregressive_model}
in a
nonparametric setting based on H\"{o}lder smoothness assumptions on
the TVAR
coefficients. A~lower bound of the prediction risk is given in this
setting and
this result is used to show that the proposed aggregation methods
achieve the
minimax adaptive rate. Section~\ref{subsection:oracle_inequalities} contains
the proofs of the oracle inequalities. The proof of the lower bound of
the minimax
prediction risk is presented in
Section~\ref{section:proof_lower_bound}.
Numerical experiments illustrating these results are then described in
Section~\ref{sec:numer-exper}.
One \hyperref[section:useful_bounds_TVAR_postponed_proofs]{Appendix} and one supplementary material [\citet
{GiraudRoueffSanchez-Perez2015s}] complete this
paper. \hyperref[section:useful_bounds_TVAR_postponed_proofs]{Appendix} and
[\citeauthor{GiraudRoueffSanchez-Perez2015s} (\citeyear{GiraudRoueffSanchez-Perez2015s}), Section~A]
contain some
postponed proofs and useful lemmas, [Giraud, Roueff and Sanchez-Perez 
(\citeyear{GiraudRoueffSanchez-Perez2015s}), Section~B] explains how
to build nonadaptive minimax
predictors which can be used in the aggregation step and
[Giraud, Roueff and Sanchez-Perez (\citeyear
{GiraudRoueffSanchez-Perez2015s}), Section~C] provides additional
results with
improved aggregation rates.

\section{Online aggregation of predictors for nonstationary processes}
\label{sec:online-aggr-pred}

\subsection{General model}

In this section, we consider a time series $(X_t)_{t\in\Z}$ admitting
the following
\emph{nonstationary} sub-linear property with respect to the nonnegative
process $(Z_t)_{t\in\Z}$.

\begin{longlist}[(M-1)]\label{hyp:gen}
\item[(M-1)] The process $(X_{t})_{t\in\Z}$ satisfies
%
\begin{equation}
\label{eq:gen-series-representation} \llvert X_{t}\rrvert\leq\sum
_{j\in\Z}A_{t}(j)  Z_{t-j}  ,
\end{equation}
where $ (A_{t}(j) )_{t,j\in\Z}$ are nonnegative coefficients
such that
%
\begin{equation}
\label{eq:unif-coeff-bound} A_*: = \sup_{t\in\Z}\sum
_{j\in\Z} A_t(j) <\infty .
\end{equation}
\end{longlist}
Additional assumptions will be required on $(Z_t)_{t\in\Z}$ to deduce
useful properties for $(X_t)_{t\in\Z}$. Note, for instance, that the
condition on
$A_*$ in~(\ref{eq:unif-coeff-bound}) guarantees that, if $(Z_t)_{t\in
\Z}$ has a
uniformly bounded $L^p$-norm, the convergence of the infinite sum
in~(\ref{eq:gen-series-representation}) holds almost surely and in the
$L^p$-sense, with both convergences defining the same limit. It follows that
$(X_t)_{t\in\Z}$ also has uniformly bounded $L^p$ moments. Let us
give some
particular contexts where the representation~(M-1) can be
used.

\begin{example}[(Time varying linear processes)]\label{exple:tvlinear}
Standard weakly stationary processes such as ARMA processes [see
\citet{BrockwellDavis2006}] admit a Wold decomposition of the form
\[
X_t = \sum_{j\geq0} a(j)
\xi_{t-j}  , %
\]
where $(\xi_{t})_{t\in\Z}$ is a weak white noise with, says, unit
variance. This model, sometimes referred to as an MA($\infty$)
representation, is often extended to a two-sided sum representation
\[
X_t = \sum_{j\in\Z} a(j)
\xi_{t-j}  , %
\]
and additional assumptions on the existence of higher moments for
$(\xi_{t})_{t\in\Z}$ or on the independence of the $\xi_t$'s are
often used for
statistical inference or prediction; see
\citet{BrockwellDavis2006}, Chapters~7 and~8.
Because
the sequence $(A_t(j))_{j\in\Z}$ may vary with $t$ in~(M-1),
we may
extend this standard stationary setting and also consider
linear processes with time varying coefficients. In this
case, we have
%
\begin{equation}
\label{eq:tvlin-rep} X_{t} = \sum_{j\in\Z}a_{t}(j)
 \xi_{t-j}  ,
\end{equation}
where $(\xi_t)$ is a sequence of centered independent random variables with
unit variance and $(a_{t}(j))_{t,j\in\Z}$ is supposed to
satisfy~(\ref{eq:unif-coeff-bound}) with $A_t(j)=\llvert
a_{t}(j)\rrvert $, so
that~\mbox{(M-1)} holds with $Z_t=\llvert \xi_t\rrvert $.
For this general class of processes, statistical inference is not
easily carried out: each new observation $X_t$ comes with a new
unknown sequence $(a_t(j))_{j\in\Z}$. However, additional assumptions
on this
set of sequences allow to derive and study appropriate statistical inference
procedures. A sensible approach in this direction is to consider a
\emph{locally
stationary} model as introduced in \citet{Dahlhaus1996}. In this
framework, the
set of sequences
$\{(a_t(j))_{j\in\Z}, 1\leq t\leq T\}$ is controlled as $T\to\infty
$ by
artificially (but meaningfully) introducing a dependence in $T$, hence
is written as
$(a_{t,T}(j))_{j\in\Z,1\leq t\leq T}$, and by approximating it with a
set of sequences rescaled on the time interval $[0,1]$, $a(u,j)$, $u\in[0,1]$,
$j\in\Z$, for example, in the following way:
\[
\sup_{T\geq1}\sup_{j\in\Z}\sum
_{t=1}^T\biggl\llvert a_{t,T}(j)-a
\biggl(\frac{t}{T},j \biggr)\biggr\rrvert< \infty . %
\]
Then various interesting statistical inference problems based on
$X_1,\dots,X_T$ can be tackled by assuming some smoothness on the mapping
$u\mapsto a(u,j)$ and, possibly, additional assumptions on the
structure of the
sequence\break $(a(u,j))_{j\in\Z}$ for each $u\in[0,1]$ [see \citet
{Dahlhaus2009} and
the references therein].
\end{example}

%
\begin{example}[(TVAR model)]\label{exple:tvar}
A particular instance of Example~\ref{exple:tvlinear} is the
\emph{time varying autoregressive} (TVAR) process, which is assumed to satisfy
the recursive equation
\[
X_t = \sum_{j=1}^d
\theta_{j,t} X_{t-j} + \sigma_t
\xi_t  , %
\]
where $(\xi_t)_{t\in\Z}$ is a white noise process; see \citet
{Grenier1983}. It
turns out that, in the framework introduced by \citet
{Dahlhaus1996}, under
suitable assumptions, such processes admit a time varying linear representation
of the form~(\ref{eq:tvlin-rep}); see \citet
{Kunsch1995,Dahlhaus1996}. In
Section~\ref{section:time-varying_autoregressive_model}, we focus on
such a
class of processes and use the aggregation of predictors to derive adaptive
minimax predictors under specific smoothness assumptions on the time varying
coefficients.
\end{example}

%
\begin{example}[(A nonlinear extension)] \label{example:non-linear}
It can also be interesting to consider nonlinear extensions of
Example~\ref{exple:tvar}. A simple example is obtained by setting
\[
X_t=g_{t}(X_{t-1})+\xi_t  ,
\]
where $(\xi_t)_{t\in\Z}$ is an i.i.d. sequence and $g_t$
is a time varying sub-linear sequence of functions satisfying, for all
$t$ that
\[
\bigl\llvert g_{t}(x)\bigr\rrvert\leq\alpha\llvert x\rrvert ,
\]
for some $\alpha\in(0,1)$. Since $g_t$ is no longer linear but
sub-linear, such
a model does not enjoy an exact linear representation of the
form~(\ref{eq:tvlin-rep}). Nevertheless, since we have
\[
\llvert X_t\rrvert\leq\alpha\llvert X_{t-1}\rrvert+
\llvert\xi_t\rrvert , %
\]
and\vspace*{1pt} iterating this equation backwards yields assumption~(M-1) with
$Z_t=\llvert \xi_t\rrvert $ and $A_t(j)=\alpha^j$. In the stationary
case, where
$g=g_t$ does
not depend on $t$, a well-known nonlinear extension is the threshold
autoregressive model where $g$ is piecewise linear; see
\citet{TongLim1980}.
\end{example}

Our goal in this section is to derive oracle bounds for the aggregation of
predictors that hold for the general model~(M-1) with one of
the two
following additional assumptions on $(Z_t)_{t\in\Z}$.
\begin{longlist}[(N-1)]
\item[(N-1)] The nonnegative process $(Z_{t})_{t\in\Z}$
satisfies
\[
\m_p:=\sup_{t\in\Z}\E\bigl[Z_{t}^p
\bigr]<\infty . %
\]
\item[(N-2)] The nonnegative process
$(Z_{t})_{t\in\Z}$ is a sequence of independent random variables
fulfilling
\[
\mphi(\zeta):=\sup_{t\in\Z}\E\bigl[\rme^{\zeta Z_{t}} \bigr]<
\infty . %
\]
\end{longlist}
Assumptions~\textup{(N-1)}
and~\textup{(N-2)} appear to be quite mild. As
mentioned in
Example~\ref{exple:tvlinear}, basic assumptions in stationary time
series usually
include moments of sufficiently high order for the innovations and their
independence, or rely on the Gaussian assumption, which is contained
in~\textup{(N-2)}. We also note that, in the context of
locally stationary time series, our assumptions on the innovations are weaker
than those used in the recent works
\citet{Dahlhaus2009}, \citeauthor{DahlhausPolonik2006} (\citeyear{DahlhausPolonik2006,DahlhausPolonik2009}). Precise
comparisons between our assumptions and usual ones in the aggregation
literature will be given after
Corollary~\ref{corollary:oracle-bounds}.

\subsection{Aggregation of predictors} \label{subsection:aggregation}
Let $(x_{t})_{t\in\Z}$ be a real valued sequence. We say that
$\widehat{x}_{t}$ is a predictor of $x_{t}$ if it is a measurable
function of
$(x_s)_{s\leq t-1}$. Throughout this paper, the quality of a sequence of
predictors
$(\widehat{x}_{t})_{1\leq t\leq T}$ is evaluated for some  $T\geq1$
using the
$\ell^2$ loss averaged over the time period $\{1,\dots,T\}$
\[
\frac{1}T\sum_{t=1}^T (
\widehat{x}_t-x_t )^2. %
\]
Now, given a
collection of $N$ sequences of predictors $ \{(\widehat
{x}_{t}^{(i)})_{1\leq t\leq T},
1\leq i\leq N \}$,
we wish to sequentially derive a new predictor which predicts almost
as accurately as or more accurately than the best of them.

\begin{algorithm}[b]
\caption{Online computation of the aggregation algorithms}
\label{algo:strat12}
\KwParameters{the learning rate $\eta$ (in $(0,\infty)$) and the
strategy (1 or 2)}
\KwInit{$t=1$, $\widehat{\alpha}_{t}=(1/N)_{i=1,\dots,N}$}
\While{\KwIn{the predictions $\widehat{x}^{(i)}_{t}$ for
$i=1,\dots,N$}}{
$\widehat{x}_{t}=\widehat{x}_{t}^{[\widehat{\alpha}_{t}]}=\sum
_{i=1}^{N}\widehat{\alpha}_{i,t}\widehat{x}_{t}^{(i)}$;\\
\Return{$\widehat{x}_{t}$;}\\

\KwAndWhen{\KwIn{a new $x_{t}$}}
{
$t=t+1$;\\
\For{$i= 1$ \KwTo$N$}{
\Switch{strategy}{
{$\mathbf{case}$ 1: $v_{i,t}=\widehat{\alpha}_{i,t-1}\exp(-2\eta
(\widehat{x}^{[\widehat{\alpha
}_{t-1}]}_{t-1}-x_{t-1} )\widehat{x}^{(i)}_{t-1} )$ }\\
{$\mathbf{case}$ 2: $v_{i,t}=\widehat{\alpha}_{i,t-1}\exp(-\eta
(\widehat{x}^{(i)}_{t-1}-x_{t-1} )^2 )$ }
}
$\widehat{\alpha}_{t}= (v_{i,t}/\sum_{k=1}^{N}v_{k,t} )_{i=1,\ldots,N}$;
}
}
}
\end{algorithm}

In the present paper and for our purposes, aggregating the predictors
amounts to compute a convex combination of them at each
time $t$. This corresponds to choosing at each time $t$ an element
$\alpha_{t}$ of the simplex
%
%
\begin{equation}
\label{eq:def-simplexN} \mathcal{S}_{N}= \Biggl\{{\mathbf s}=(s_{1},
\ldots,s_{N})\in\R_{+}^{N}:\sum
 _{i=1}^{N}s_{i}=1 \Biggr\}
\end{equation}
and compute
\[
\widehat{x}_t^{[\alpha_{t}]}=\sum_{i=1}^N
\alpha_{i,t} \widehat{x}_t^{(i)} . %
\]
We consider two strategies of aggregation, which are studied in the
context of
bounded sequences in \citet{Cesa-BianchiLugosi2006,Catoni2004}.
More recent
contributions and extensions can be found in \citet{Gerchinovitz2011}.
See also
\citet{Stoltz2011} for a pedagogical \hyperref[sec1]{Introduction}. These
strategies are
sequential and online, meaning that:
\begin{longlist}[(ii)]
\item[(i)] to compute the aggregation
weights $\alpha_{t}$ at time $t$, only the values of $\{\widehat
{x}_{s}^{(i)}, 1\leq i\leq
N\}$ and $x_s$ up to time $s=t-1$ are used,

\item[(ii)] the computation can be done
recursively by updating a set of quantities, the number of which does
not depend on $t$.
\end{longlist}

These two properties are met in the Algorithm~\ref{algo:strat12}
detailed below.

We consider in the remaining of the paper a convex aggregation of predictors
\[
\widehat{x}_{t} =\widehat{x}_{t}^{[\widehat{\alpha}_{t}]}= \sum
 _{i=1}^{N}\widehat{\alpha}_{i,t}\widehat
{x}_{t}^{(i)},\qquad1\leq t\leq T , %
\]
with some specific weights $\widehat{\alpha}_{i,t}$ defined as follows.

\subsubsection*{Strategy 1: Building weights from the gradient of the
quadratic loss}
The first strategy is to define for all
$i=1,\dots,N$ and $t=1,\dots,T$, the weights $\widehat{\alpha
}_{i,t}$ by
%
\begin{eqnarray}
\widehat{\alpha}_{i,t} &=& \frac{\exp
(-2\eta\sum _{s=1}^{t-1} (\sum _{j=1}^{N}\widehat{\alpha}_{j,s}\widehat{x}_{s}^{(j)}-x_{s} )\widehat
{x}_{s}^{(i)} )}{ \sum _{k=1}^{N}\exp(-2\eta\sum _{s=1}^{t-1} (\sum
_{j=1}^{N}\widehat{\alpha}_{j,s}\widehat
{x}_{s}^{(j)}-x_{s} )\widehat{x}_{s}^{(k)} )} ,
\label
{definition:alphas_gradient_quadratic_loss}
\end{eqnarray}
with the convention that a sum over no element is zero, so
$\widehat{\alpha}_{i,1}=1/N$ for all $i$.

The parameter $\eta>0$, usually called the \emph{learning rate}, will
be specified later.

\subsubsection*{Strategy 2: Building weights from the quadratic loss}
The second strategy is to define for all
$i=1,\dots,N$ and $t=1,\dots,T$, the weights $\widehat{\alpha
}_{i,t}$ by
%
\begin{eqnarray}
\widehat{\alpha}_{i,t} &=&  \frac{ \exp
(-\eta\sum _{s=1}^{t-1} (\widehat
{x}_{s}^{(i)}-x_{s} )^{2} )}{ \sum
_{k=1}^{N}\exp(-\eta\sum _{s=1}^{t-1} (\widehat
{x}_{s}^{(k)}-x_{s} )^{2} )} ,
\label
{definition:alphas_quadratic_loss}
\end{eqnarray}
with again the convention that a sum over no element is zero.

Both strategies yield the same algorithm up to
the line where $v_{i,t}$ is computed. For sake of brevity, we write
only one algorithm (see Algorithm~\ref{algo:strat12}) and use a
switch/case statement to distinguish between the two
strategies. Note, however, that the choice of the strategy ($1$ or
$2$) holds for the whole sequence of predictions.

\subsection{Oracle bounds}\label{subsection:general:bound}
We establish oracle bounds on the average prediction error of the aggregated
predictors. These bounds ensure that the error is equal to that
associated with the best convex combination of the predictors or with
the best predictor (depending on the aggregation strategy), up to two
remaining terms. One remaining term depends on the number $N$ of
predictors to aggregate and the
other one on the \emph{variability} of the original process. The
learning rate
$\eta$ can then be chosen to achieve the best trade-off between these
two terms.

The second remaining term indirectly depends on the variability of the
predictors. We control below this variability in terms of the
variability of the original process by using the
following Lipschitz property.

\begin{definition}
\label{def:unif-lip-pred}
Let $L=(L_s)_{s\geq1}$ be a sequence of nonnegative numbers. A~predictor $\widehat{x}_{t}$ of $x_t$ from $(x_s)_{s\leq t-1}$ is said
to be
$L$-Lipschitz if
\[
\llvert\widehat{x}_{t}\rrvert\leq\sum
_{s\geq1} L_s \llvert x_{t-s}\rrvert .
\]
\end{definition}

We more specifically consider a sequence $L$ satisfying the following
assumption.
\begin{longlist}
\item[(L-1)] The sequence $L=(L_s)_{s\geq1}$ satisfies
%
\begin{equation}
\label{eq:hyp:lip} L_* = \sum_{j\geq1} L_j <
\infty .
\end{equation}
\end{longlist}

This condition is trivially satisfied by constant
linear predictors depending only on a finite number of previous observations,
that is, $\widehat{x}_{t}=\sum_{s=1}^{d} L_s x_{t-s}$. In
\citeauthor{GiraudRoueffSanchez-Perez2015s} [(\citeyear
{GiraudRoueffSanchez-Perez2015s}), Section~B.1], we extend this case in the
context of the
TVAR process where the coefficients $L_s$ are replaced by estimates of the
time varying autoregressive coefficients. More generally,
assumption~\mbox{\textup{(L-1)}} appears to be quite natural in the general context
where $\E[X_t\mid(X_{t-s})_{s\geq1}]=f_t((X_{t-s})_{s\geq1})$, where
$f_t$ is a
Lipschitz function from $\R^{\N^{*}}$ to $\R$, with Lipschitz
coefficients satisfying
a condition similar to~(\ref{eq:hyp:lip}); see, for instance,
\citet{DoukhanWintenberger2008} in the case of stationary time series.

We now state two upper-bounds on the mean quadratic prediction error of the
aggregated predictors defined in the previous section, when the process $X$
fulfills the sub-linear property~(M-1).\vspace*{-1pt}

%
\begin{theorem}\label{thm:oracle-bounds}
Assume that assumption~\textup{(M-1)} holds. Let
$\{(\widehat{X}_t^{(i)})_{1\leq t\leq T}$, $ 1\leq i \leq N\}$ be a
collection of sequences of $L$-Lipschitz predictors with $L$
satisfying~\mbox{\textup{(L-1)}}.
\begin{longlist}[(iii)]
\item[(i)] Assume that the noise $Z$
fulfills~\textup{(N-1)} with $p=4$ and
let $\widehat{X}=(\widehat{X}_t)_{1\leq t\leq T}$ denote the aggregated
predictor obtained using the
weights~(\ref{definition:alphas_gradient_quadratic_loss}) with any
$\eta>0$. Then we have\vspace*{-1pt}
\begin{eqnarray}
\label{eq:strat1-final} {1\over T}\sum_{t=1}^T
\E\bigl[ (\widehat X_{t}-X_{t} )^2 \bigr]
&\leq&\inf_{\nu\in\mathcal{S}_N} {1\over T}\sum
_{t=1}^T\E\bigl[ \bigl(\widehat
X_t^{[\nu]}-X_{t} \bigr)^2 \bigr]
\nonumber\\[-9pt]\\[-9pt]\nonumber
&&{} +{\log N\over T\eta}+2\eta (1+L_*)^{4} A_{*}^4 \m_{4} .
\end{eqnarray}

\item[(ii)]  Assume that the noise $Z$
satisfies~\textup{(N-1)} with a given $p>2$
and
let $\widehat{X}=(\widehat{X}_t)_{1\leq t\leq T}$ denote the aggregated
predictor obtained using the
weights (\ref{definition:alphas_quadratic_loss}) with any $\eta>0$.
Then we have
\begin{eqnarray}
\label{eq:strat2_p-final} {1\over T}\sum_{t=1}^T
\E\bigl[ (\widehat X_{t}-X_{t} )^2 \bigr]
&\leq& \min_{1\leq i\leq N} {1\over T}\sum
_{t=1}^T\E\bigl[ \bigl(\widehat
X_{t}^{(i)}-X_{t} \bigr)^2 \bigr]
\nonumber\\[-10pt]\\[-10pt]\nonumber
&&{} + \frac{\log N}{T\eta} + (2\eta)^{p/2-1}A_{*}^{p}(1+L_*)^{p}m_{p} .
\end{eqnarray}

\item[(iii)] Assume that the noise $Z$
fulfills~\textup{(N-2)} for some positive $\zeta$
and
let $\widehat{X}=(\widehat{X}_t)_{1\leq t\leq T}$ denote the aggregated
predictor obtained using the
weights~(\ref{definition:alphas_quadratic_loss}) with $\eta>0$.
Then, for any
%
\begin{equation}
\label{eq:def-sup-unif-coeff-bound} \lambda\in\biggl(0,\frac{\zeta
}{a^*(L_*+1)} \biggr]\qquad\mbox{with } a^*: = \sup
_{j\in\Z}\sup_{t\in\Z} A_t(j)  \leq
  A_{*} ,
\end{equation}
we have
\begin{eqnarray}\label{eq:strat2_exp-final}
&& {1\over T}\sum_{t=1}^T
\E\bigl[ (\widehat X_{t}-X_{t} )^2 \bigr]\nonumber
\\
&&\qquad \leq\min_{1\leq i\leq N} {1\over T}\sum
_{t=1}^T\E\bigl[ \bigl(\widehat
X_{t}^{(i)}-X_{t} \bigr)^2 \bigr]
\\
&&\quad\qquad{} + \frac{\log N}{T\eta} + \frac{2}{\rme}\lambda^{-2} \bigl(2+\lambda(2
\eta)^{-1/2} \bigr)  \rme^{- \lambda(2\eta)^{-1/2}}  \bigl(\mphi(\zeta)
\bigr)^{\lambda A_*(1+L_*)/\zeta} .\nonumber
\end{eqnarray}
\end{longlist}
\end{theorem}

The proof can be found in Section~\ref{sec:proof-theor-strat12}.

%
\begin{remark}
\label{rem:oracle-bounds-explicites}
The bounds~(\ref{eq:strat1-final}),~(\ref{eq:strat2_p-final}) and
(\ref{eq:strat2_exp-final}) are explicit in
the sense that all the constants appearing in them are directly derived from
those appearing in assumptions~(M-1), \textup{(L-1)}, \textup{(N-1)} and~\textup{(N-2)}.
\end{remark}

The following corollary is obtained by choosing $\eta$ [and
$\lambda$ in the case~(iii)] adequately in the three
cases of Theorem~\ref{thm:oracle-bounds}.

%
\begin{corollary} \label{corollary:oracle-bounds}
Assume that assumption~\textup{(M-1)} holds. Let
$\{(\widehat{X}_t^{(i)})_{1\leq t\leq T}$, $ 1\leq i \leq N\}$ be a
collection of sequences of $L$-Lipschitz predictors with $L$
satisfying~\mbox{\textup{(L-1)}}.
\begin{longlist}[(iii)]
\item[(i)] Assume that the noise $Z$
fulfills~\textup{(N-1)} with $p=4$ and
let $\widehat{X}=(\widehat{X}_t)_{1\leq t\leq T}$ denote the aggregated
predictor obtained using the
weights~(\ref{definition:alphas_gradient_quadratic_loss}) with
%
\begin{equation}
\label{eq:strat1-choice-eta} \eta= \frac{1}{ (2m_4 )^{1/2}
(1+L_*)^{2} A_{*}^2}  \biggl(\frac{\log
N}{T}
\biggr)^{1/2}  .
\end{equation}
This gives
%
\begin{equation}
\label{eq:strat1-final-eta-opt}\qquad {1\over T}\sum_{t=1}^T
\E\bigl[ (\widehat X_{t}-X_{t} )^2 \bigr]
\leq\inf_{\nu\in\mathcal{S}_N} {1\over T}\sum
_{t=1}^T\E\bigl[ \bigl(\widehat
X_t^{[\nu]}-X_{t} \bigr)^2 \bigr] +
C_1  \biggl(\frac{\log
N}{T} \biggr)^{1/2} ,
\end{equation}
with $C_1=2 (2m_4 )^{1/2} (1+L_*)^{2} A_{*}^2$.

\item[(ii)]  Assume that the noise $Z$
satisfies~\textup{(N-1)} with a given $p> 2$
and
let $\widehat{X}=(\widehat{X}_t)_{1\leq t\leq T}$ denote the aggregated
predictor obtained using the
weights~(\ref{definition:alphas_quadratic_loss}) with
%
\begin{equation}
\label{eq:strat2_p-choice-eta} 
\eta= \frac{1}{2m_{p}^{2/p}(1+L_*)^{2} A_{*}^{2}}  \biggl(
\frac
{\log N}{T} \biggr)^{2/p}.
\end{equation}
We then have
%
\begin{eqnarray}\label{eq:strat2_p-final-eta-opt}
&& {1\over T}\sum_{t=1}^T
\E\bigl[ (\widehat X_{t}-X_{t} )^2 \bigr]
\nonumber\\[-8pt]\\[-8pt]\nonumber
&&\qquad \leq\min_{1\leq i\leq N} {1\over T}\sum
_{t=1}^T\E\bigl[ \bigl(\widehat
X_{t}^{(i)}-X_{t} \bigr)^2 \bigr] +
C_2  \biggl(\frac{\log
N}{T} \biggr)^{1-2/p} ,
\end{eqnarray}
with $C_2=3m_{p}^{2/p}(1+L_*)^{2} A_{*}^{2}$.

\item[(iii)] Assume that the noise $Z$
fulfills~\textup{(N-2)} for some positive $\zeta$
and
let $\widehat{X}=(\widehat{X}_t)_{1\leq t\leq T}$ denote the aggregated
predictor obtained using the
weights~(\ref{definition:alphas_quadratic_loss}) with
%
\begin{equation}
\label{eq:strat2_exp-choice-eta} \eta= \frac{\zeta
^2}{2(1+L_*)^2A_{*}^2}  \biggl(\log\biggl(\frac
{T}{\log N}
\biggr) \biggr)^{-2} .
\end{equation}
Then we have
\begin{eqnarray}
\label{eq:strat2_exp-final-opt-eta}
&& {1\over T}\sum_{t=1}^T
\E\bigl[ (\widehat X_{t}-X_{t} )^2 \bigr]\nonumber
\\
&&\qquad \leq\min_{1\leq i\leq N} {1\over T}\sum
_{t=1}^T\E\bigl[ \bigl(\widehat
X_{t}^{(i)}-X_{t} \bigr)^2 \bigr]+ \frac{2A_{*}^2(L_*+1)^2}{\zeta^2}  \frac{\log N}{T}
\\
&&\quad\qquad{}  \times \biggl\{
\biggl(\log\biggl(
\frac{T}{\log N} \biggr) \biggr)^{2} +\frac{\phi(\zeta)}\rme\biggl(2+
\log\biggl(\frac{T}{\log N} \biggr) \biggr) \biggr\} .\nonumber
\end{eqnarray}
[Note that when $(\log N)/T\to0$, the term between curly brackets is
equivalent to $(\log(T/\log N))^{2}$.]
\end{longlist}
\end{corollary}

Cases~(i) and~(ii) in Corollary~\ref
{corollary:oracle-bounds}
follow directly from Theorem~\ref{thm:oracle-bounds}. Case~(iii) is
more delicate since it requires optimizing $\lambda$ as well as $\eta
$ in the
second line of~(\ref{eq:strat2_exp-final}). The details are postponed to
Section~\ref{sec:proof-corollary-strat12}.

%
\begin{remark}
We\vspace*{1pt} observe that the bound in~(\ref{eq:strat2_exp-final-opt-eta}) improves
that in (\ref{eq:strat2_p-final-eta-opt}) for any $p>2$. For\vspace*{1pt} $p>4$, the
remaining term $(\log N/T)^{1-2/p}$ in~(\ref
{eq:strat2_p-final-eta-opt}) is
smaller than the remaining term $(\log N/T)^{1/2}$
in~(\ref{eq:strat1-final-eta-opt}). Similarly, the remaining term
$\log N (\log T)^2/T$ in~(\ref{eq:strat2_exp-final-opt-eta}) is
smaller than
$(\log N/T)^{1/2}$ in~(\ref{eq:strat1-final-eta-opt}). Yet, we emphasize
that the oracle inequalities~(\ref{eq:strat2_p-final-eta-opt})
and~(\ref{eq:strat2_exp-final-opt-eta}) compare the prediction risk of
$\widehat X$ to the prediction\vspace*{1pt} risk of the \emph{best predictor}
$\widehat
X^{(i)}$, while\vspace*{1pt} the oracle inequality~(\ref{eq:strat1-final-eta-opt})
compare the
prediction risk of $\widehat X$ to the prediction risk of the \emph{best
convex combination of the predictors} $\widehat X^{(i)}$, so they
cannot be
directly compared.
\end{remark}

%
\begin{remark}\label{improve-with-net}
As explained in \citeauthor{GiraudRoueffSanchez-Perez2015s} [(\citeyear
{GiraudRoueffSanchez-Perez2015s}), Section~C], under the
hypotheses of cases~(ii) and~(iii) and for certain
values of $T$~and~$N$, using a more involved aggregation step, we can get a new predictor
satisfying an oracle inequality better than that
in~(\ref{eq:strat1-final-eta-opt}). For example, under the hypotheses of
case~(iii), for $T>N^2(\log T)^6$, the remaining term
$(\log N/T)^{1/2}$ in~(\ref{eq:strat1-final-eta-opt}) can be replaced by
$N(\log T)^3/T$ which is smaller; see \citeauthor
{GiraudRoueffSanchez-Perez2015s} [(\citeyear
{GiraudRoueffSanchez-Perez2015s}), inequality~(C.7),
page~8]. Yet, this aggregation has a prohibitive
computational cost and seems difficult to implement in practice.
\end{remark}

%
\begin{remark}
\label{rem:eta-choice-VS-conditions}
In cases~(ii) and~(iii), which correspond to the
weights~(\ref{definition:alphas_quadratic_loss}), the choice of the optimal
$\eta$ depends on the assumptions on the noise,
namely \mbox{\textup{(N-1)}} or \textup{(N-2)}. Under a
moment condition of order $p$, the optimal $\eta$ is of order $(\log
N/T)^{2/p}$ and under
an exponential condition, it is of order $(\log T)^{-2}$. It is known from
\citeauthor{Catoni2004} [(\citeyear{Catoni2004}), Proposition~2.2.1]
and
\citeauthor{Yang2004} [(\citeyear{Yang2004}), Theorem~5]
that $\eta$ can be chosen as a constant (provided that it is small
enough) under a bounded noise condition, or under an exponential moment
condition on the noise for predictors at a bounded distance from the
conditional mean. Hence,
coarsely speaking, the heavier the tail of the noise, the smallest
$\eta$
should be chosen. Observing that $\eta$ allows us to tune the
influence of
the empirical risk on the weights from no influence at all ($\eta=0$ yielding
uniform weights) to the selection of the empirical risk minimizer
($\eta\to\infty$), the specific choices of $\eta$ can be
interpreted as
follows: the heavier the tail of the noise, the less we can trust the
empirical risk.
\end{remark}

\subsubsection*{Comparison with previous works} \label{section:comparison}
In the literature, prediction risk bounds of the
form~(\ref{eq:strat1-final-eta-opt}) [case~(i) of
Corollary~\ref{corollary:oracle-bounds}] are sometimes called \emph{convex
regret bounds}, and prediction risk bounds of the
form~(\ref{eq:strat2_p-final-eta-opt}) and~(\ref{eq:strat2_exp-final-opt-eta})
[cases~(ii) and~(iii) of
Corollary~\ref{corollary:oracle-bounds}] are sometimes called \emph
{best predictor
regret bounds}.

\citet{Sancetta2010} exhibits convex regret bounds in a setting close
to ours,
namely for an online aggregation of predictors for a sequence of possibly
dependent random variables. Under our moment
condition~\textup{(N-1)} with $p=4$,
\citeauthor{Sancetta2010} [(\citeyear{Sancetta2010}), Theorem~2]
provides an upper bound similar
to~(\ref{eq:strat1-final-eta-opt}) but with our remaining term $(\log
N/T)^{1/2}$ replaced by $(N\log(N)/T)^{1/2}$. Under the exponential
condition~\mbox{\textup{(N-2)}},
\citeauthor{Sancetta2010} [(\citeyear{Sancetta2010}), Theorem~1]
provides an\vspace*{1pt} upper bound similar to~(\ref{eq:strat1-final-eta-opt})
but with a remaining term $(\log N/T)^{1/2}\times(\log(NT))^2$, which
is still larger
than our remaining term under moment conditions.

Best predictor regret bounds can be found in \citet{Yang2004} for some
sequences of possibly dependent random variables.
The predictors are assumed to remain at
a bounded distance to the conditional means and the scaled
innovation noise is assumed to have either a known
distribution (satisfying a certain technical condition) or an
exponential moment. The regret
bounds are presented in a slightly different fashion from ours but it
is easy to
see that a similar result as our bound~(\ref
{eq:strat2_exp-final-opt-eta}) is
obtained in this setting. However, we do not require bounded prediction
errors and our conditions on the
noise are milder.

The i.i.d. setting has received much more attention and, even if the
setting is quite different, it is interesting to
briefly compare our results to previous works in this case. Let us
start with
the convex regret bound in case~(i) of
Corollary~\ref{corollary:oracle-bounds}. Most of the existing results
[see, e.g., \citet{JuditskyNemirovski2000}, \citet{Yang2000}, \citet{Tsybakov2003} or
\citet{WangPaterliniGaoYang2014} for recent extensions to $\ell^q$
aggregation] assume the predictors to be bounded and various conditions
on the noise are
considered (very often the noise is assumed to be Gaussian).
In such settings, the best possible remaining term typically takes the form
$(\log N/T)^{1/2}$ when $N$ is much larger than $T^{1/2}$ and of the
form $N/T$
if $N$ is smaller than $T^{1/2}$; see \citeauthor
{JuditskyNemirovski2000} [(\citeyear{JuditskyNemirovski2000}),
Theorem~3.1],
\citeauthor{Yang2004} [(\citeyear{Yang2004}), Theorem~6]
and \citeauthor{Tsybakov2003} [(\citeyear{Tsybakov2003}), Theorem~2].
Hence, our
bound~(\ref{eq:strat1-final-eta-opt}) is similar only in the case
where $N$ is
much larger than $T^{1/2}$.
However, as explained in Remark~\ref{improve-with-net} and
[Giraud, Roueff and Sanchez-Perez (\citeyear
{GiraudRoueffSanchez-Perez2015s}), Section~C], when $T$ is larger than
$N^2$ and under
the moment
condition~\textup{(N-2)}, we can get via a more
involved aggregation procedure, a convex regret bound with a remaining
term of the same order $N/T$ up to a $(\log T)^3$ factor
[see \citet{GiraudRoueffSanchez-Perez2015s},
inequality (C.7), page~8]. Let us now compare our
bound~(\ref{eq:strat2_p-final-eta-opt}) in case~(ii) to
optimal bounds in the i.i.d. setting under moment conditions on the
noise. Corollary~7.2 and Theorem~8.6 in~\citet{Audibert2009} shows
that the optimal aggregation rate is $(\log N/T)^{1-2/(p+2)}$ in the
i.i.d. setting with bounded predictors and moment conditions of order
$p$ on the noise.
Our remaining term $(\log N/T)^{1-2/p}$ in~(\ref
{eq:strat2_p-final-eta-opt}) is slightly larger, yet an inspection of
the proof of \citeauthor{Audibert2009} [(\citeyear{Audibert2009}),
Corollary~7.2] shows that the
aggregation rate would also be $(\log N/T)^{1-2/p}$ in this corollary,
if the predictors were assumed to have a moment condition of order $p$
instead of being uniformly bounded (we are not aware of any lower bound
in this setting matching this rate). Finally, when the data and the
predictors are bounded, the best aggregation rate is known to be $(\log
N)/T$ in the i.i.d. setting; see, for example,
\citet{Audibert2009}, Theorem~8.4.
Our bound~(\ref{eq:strat2_exp-final-opt-eta}) in case~(iii)
achieves the same rate up to a $(\log T)^2$ factor.

\section{Time-varying autoregressive (TVAR) model} \label
{section:time-varying_autoregressive_model}

\subsection{Nonparametric TVAR model}
\subsubsection{Vector norms and H\"older smoothness norms}
We introduce some preliminary notation before defining the model.
In the remainder of this\vadjust{\goodbreak} article, vectors are denoted using boldface
symbols and
$\llvert \mathbf{x}\rrvert $ denotes the Euclidean norm of $\mathbf{x}$,
$\llvert \mathbf{x}\rrvert =(\sum_i\llvert x_i\rrvert ^2)^{1/2}$.

For $\beta\in(0,1]$ and an interval $I\subseteq\R$, the $\beta$-H\"{o}lder
semi-norm of a function $\mathbf{f}:I\rightarrow\R^{d}$ is defined by
\[
\llvert\mathbf{f}\rrvert_{\beta} = \sup
_{0<\llvert s-s'\rrvert <1}{
 \frac
{\llvert \mathbf{f}(s)-\mathbf{f}(s')\rrvert }{\llvert s-s'\rrvert
^{\beta}}}  . %
\]
This semi-norm is extended to any $\beta>0$ as follows. Let $k\in\N$
and $\alpha\in(0,1]$ be such that $\beta=k+\alpha$. If
$\mathbf{f}$ is $k$ times differentiable on $I$, we define
\[
\llvert\mathbf{f}\rrvert_{\beta} = \bigl\llvert\mathbf{f}^{(k)}
\bigr\rrvert_{\alpha}  , %
\]
and $\llvert \mathbf{f}\rrvert _{\beta}=\infty$ otherwise. We
consider the case
$I=(-\infty,1]$. For $R>0$ and
$\beta>0$, the $(\beta,R)$-H\"{o}lder ball is denoted by
\[
\Lambda_{d}(\beta,R) = \bigl\{\mathbf{f}:(-\infty,1]\rightarrow\R
^{d}, \mbox{ such that } \llvert\mathbf{f}\rrvert_{\beta}
\leq R \bigr\}  . %
\]

\subsubsection{TVAR parameters in rescaled time}

The idea of using a rescaled time with the sample size $T$ for the TVAR
parameters goes back to~\citet{Dahlhaus1996}. Since then, it has
always been
a central example of locally stationary linear processes. In this
setting, the
time varying autoregressive coefficients and variance which generate the
observations $X_{t,T}$ for $1\leq t\leq T$ are represented by functions from
$[0,1]$ to $\R^{d}$ and from $[0,1]$ to $\R_+$, respectively. The definition
sets of these functions are extended to $(-\infty,1]$ in the following
definition.

%
\begin{definition}[(TVAR model)]
\label{definition:tvar}
Let $d\geq1$. Let $\theta_1,\dots,\theta_d$ and $\sigma$ be
functions defined on
$(-\infty,1]$ and $(\xi_{t})_{t\in\Z}$ be a sequence of i.i.d. random
variables with zero mean and unit variance. For any $T\geq1$, we say that
$(X_{t,T})_{t\leq T}$ is a TVAR process with time varying parameters
$\theta_1,\dots,\theta_d,\sigma^2$ sampled at frequency $T^{-1}$ and
normalized innovations $(\xi_t)$ if the two following assertions hold:
\begin{longlist}[(ii)]
\item[(i)]  The process $X$ fulfills the time varying
autoregressive equation
%
\begin{eqnarray}
\label{equation:definition_TVAR} X_{t,T} &=& \sum
_{j=1}^{d}
\theta_{j} \biggl({ \frac{t-1}{T}}
\biggr)X_{t-j,T}+\sigma\biggl({ \frac
{t}{T}} \biggr)
\xi_{t}\qquad\mbox{for } -\infty<t\leq T .
\end{eqnarray}
\item[(ii)] The sequence $(X_{t,T})_{t\leq T}$ is bounded in
probability,
\[
\lim_{M\to\infty}\sup_{-\infty<t\leq T}\PP\bigl(\llvert
X_{t,T}\rrvert>M\bigr) = 0 . %
\]
\end{longlist}
\end{definition}

This definition extends the usual definition of TVAR processes, where
the time varying parameters
$\theta_1,\dots,\theta_d$ and $\sigma^2$ are assumed to be constant on
$\R_-$; see, for example, \citeauthor{Dahlhaus1996} [(\citeyear
{Dahlhaus1996}), page~144]. The TVAR
model is generally
used for the sample $(X_{t,T})_{1\leq t\leq T}$. The definition of the
process for negative times $t$ can be seen as a way to define initial
conditions for $X_{1-d,T},\dots,X_{0,T}$, which are then sufficient to compute
$(X_{t,T})_{1\leq t\leq T}$ by iterating~(\ref{equation:definition_TVAR}).
However, in the context of prediction, it can be useful to consider predictors
$\widehat{X}_{t,T}$ which may rely on historical data $X_{s,T}$
arbitrarily far
away in the past, that is, with $s$ tending to $-\infty$. To cope with
this situation,
our definition of the TVAR process $(X_{t,T})$ holds for all time indices
$-\infty<t\leq T$ and we use the following definition for predictors.

\begin{definition}[(Predictor)] 
\label{def:predictor}
For all $1\leq t\leq T$, we say that $\widehat{X}_{t,T}$ is a
predictor of
$X_{t,T}$ if it is $\mathcal{F}_{t-1,T}$-measurable, where
%
\begin{equation}
\label{eq:filtration} \mathcal{F}_{t,T} = \sigma(X_{s,T},
s=t,t-1,t-2,\dots)
\end{equation}
is the $\sigma$-field generated by
$(X_{s,T})_{s\leq t}$.
For any $T\geq1$, we denote by $\mathcal{P}_T$ the set of sequences
$\widehat{X}_T=(\widehat{X}_{t,T})_{1\leq t\leq
T}$ of predictors for $(X_{t,T})_{1\leq t\leq T}$, that is, the set of all
processes $\widehat{X}_T=(\widehat{X}_{t,T})_{1\leq t\leq
T}$ adapted to the filtration $(\mathcal{F}_{t-1,T})_{1\leq t\leq T}$.
\end{definition}

In this general framework, the time $t=1$ corresponds to the beginning of
the aggregation procedure. Such a framework applies in two practical
situations. In
the first one, we start collecting data $X_{t,T}$ at $t\geq1$ and
compute several predictors $\widehat{X}_{t,T}^{(j)}$, $j=1,\ldots,N$
from them. Thus, the resulting aggregated predictor only depends on
$(X_{s,T})_{1\leq s\leq t-1}$. A somewhat different
situation is when historical data is available beforehand the
aggregation step, so that a given predictor $\widehat{X}_{t,T}^{(j)}$
is allowed to
depend also on data $X_{s,T}$ with $s\leq0$, while the aggregation
step only
starts at $t\geq1$, and thus depends on the data $(X_{s,T})_{s\leq0}$ only
through the predictors. It is important to note that, in
contrast to the usual stationary situation, having observed the process
$X_{s,T}$ for infinitely many $s$'s in the past (for all $s\leq t-1$)
is not
so decisive for deriving a predictor of $X_{t,T}$, since observations
far away
in the past may have a completely different statistical
behavior.

\subsubsection{Stability conditions}
The next proposition proves that under standard stability conditions on
the time varying parameters
$\theta_1,\dots,\theta_d$ and $\sigma^2$, condition~(ii) in
Definition~\ref{definition:tvar} ensures the existence and uniqueness
of the
solution of equation~(\ref{equation:definition_TVAR}) for $t\leq0$
(and thus for
all $t\leq T$). We define the time varying autoregressive polynomial by
\begin{eqnarray*}
\btheta(z;u) &=& 1-\sum
_{j=1}^{d}
\theta_{j}(u)z^{j}  .
\end{eqnarray*}
Let us denote, for any $\delta>0$,
%
\begin{equation}
\label{eq:SdDelta} s_{d}(\delta)= \bigl\{\btheta:(-\infty,1]\rightarrow
\R^{d},\btheta(z;u)\neq0, \forall\llvert z\rrvert<
\delta^{-1}, u\in[0,1] \bigr\} .
\end{equation}

Define, for $\beta>0$, $R>0$, $\delta\in(0,1)$, $\rho\in[0,1]$ and
$\sigma_+>0$,
the class of parameters
\begin{eqnarray*}
&& \mathcal{C} (\beta,R,\delta,\rho,\sigma_{+} )
\\
&&\qquad = \bigl\{ (
\bolds{\theta},\sigma):(-\infty,1]\to\R^d\times[\rho\sigma_+,
\sigma_+]:\bolds{\theta}\in\Lambda_d(\beta,R)\cap
s_d(\delta) \bigr\} .
\end{eqnarray*}

The definition of the class $\mathcal{C}$ is very similar to that of
\citet{MoulinesPriouretRoueff2005}. The domain of definition in
their case is
$[0,1]$ whereas it is $(-\infty,1]$ in ours. We have the following
stability result.

%
\begin{proposition} \label{proposition:stationary_TVAR} Assume that
the time
varying AR coefficients $\theta_{1},\dots,\theta_d$ are uniformly continuous
on $(-\infty,1]$ and the time varying variance
$\sigma^{2}$ is bounded on $(-\infty,1]$. Assume
moreover that there exists $\delta\in(0,1)$ such that $\btheta\in
s_{d} (\delta)$. Then there exists $T_0\geq1$ such that,
for all
$T\geq T_0$, there exists a unique process
$(X_{t,T})_{t\leq T}$ which satisfies~\textup{(i)} and~\textup{(ii)} in
Definition~\ref{definition:tvar}. This solution admits the linear
representation
%
\begin{equation}
\label{eq:inf-series-representation} X_{t,T} = \sum_{j=0}^{\infty}a_{t,T}(j)
  \sigma\biggl({ \frac{t-j}{T}} \biggr) \xi_{t-j},
\qquad-\infty< t\leq T  ,
\end{equation}
where the coefficients $(a_{t,T}(j))_{t\leq T,j\geq0}$ satisfy that
for any
$\delta_{1}\in(\delta,1)$,
\begin{eqnarray*}
\bar{K} &=& \sup_{T\geq T_0} \sup_{-\infty<t\leq T} \sup
_{j\geq
0}  \delta_{1}^{-j}\bigl\llvert
a_{t,T}(j)\bigr\rrvert< \infty .
\end{eqnarray*}
Moreover, if
$ (\btheta,\sigma)\in\mathcal{C} (\beta,R,\delta,0,\sigma_{+} )$
for some positive constants $\beta$, $R$ and $\sigma_+$, then
the constants $T_0$ and $\bar K$ can be chosen only depending on
$\delta_{1}$, $\delta$, $\beta$
and $R$.
\end{proposition}

A proof of Proposition~\ref{proposition:stationary_TVAR} is provided in
\hyperref[section:useful_bounds_TVAR_postponed_proofs]{Appendix}.
This kind of result is classical under various smoothness assumptions
on the
parameters and initial conditions for $X_{1-k,T}$, $k=1,\dots,d$. For instance,
in \citet{DahlhausPolonik2009}, bounded variations and a constant
$\btheta$
for negative times are used for the smoothness assumption on $\btheta$
and for
defining the initial conditions. The linear
representation~(\ref{eq:inf-series-representation}) of TVAR processes
was first
obtained in the seminal papers \citet{Kunsch1995,Dahlhaus1996}. We
note that
an important consequence of Proposition~\ref
{proposition:stationary_TVAR} is
that for any $T\geq T_0$, the process $(X_{t,T})_{t\leq T}$ satisfies
assumption~(M-1) with $Z_t=\llvert \xi_t\rrvert $ and
$A_t(j)=\llvert a_{t,T}(j)
\sigma((t-j)/T )\rrvert $ for $j\geq0$. Moreover, the
constant $A_*$
in~(\ref{eq:unif-coeff-bound}) is bounded independently of $T$, and we have,
for all
$ (\btheta,\sigma)\in\mathcal{C} (\beta,R,\delta,0,\sigma_{+} )$,
%
\begin{equation}
\label{eq:boundA_*tvar} A_*\leq\frac{\bar K\sigma_+}{1-\delta_{1}}  ,
\end{equation}
where $\bar K>0$ and $\delta_{1}\in(0,1)$ can be chosen only
depending on $\delta$,
$\beta$ and $R$.

\subsubsection{Main assumptions}

Based on Proposition~\ref{proposition:stationary_TVAR}, given an
i.i.d. sequence
$(\xi_t)_{t\in\Z}$ and constants $\delta\in(0,1)$, $\rho\in[0,1]$,
$\sigma_+>0$, $\beta>0$ and $R>0$, we consider the following assumption.

\begin{longlist}[(M-2)] 
\item[(M-2)] The sequence $(X_{t,T})_{t\leq T}$ is a TVAR
process with time varying standard deviation $\sigma$, time varying AR
coefficients
$\theta_1,\dots,\theta_d$ and innovations $(\xi_{t})_{t\in\Z}$, and
$(\btheta,\sigma)\in\mathcal{C} (\beta,R,\delta,\rho,\sigma
_{+} )$.
\end{longlist}
Let $\xi$ denote a generic random variable with the same distribution
as the
$\xi_{t}$'s.
Under assumption~(M-2), the distribution of
$ (X_{t,T} )_{1-d\leq t\leq T}$ only depends on that of $\xi
$ and on
the functions $\btheta$ and $\sigma$.
For a given distribution $\psi$ on $\R$ for $\xi$, we denote by
$\PP^\psi_{(\btheta,\sigma)}$ the probability distribution of the
whole sequence
$(X_{t,T})_{t\leq T}$ and by $\E^\psi_{(\btheta,\sigma)}$ its corresponding
expectation.

The next two assumptions on the innovations are useful to prove upper
bounds of the prediction error.
\begin{longlist}[(I-3)]\label{hyp:innov-moment}
\item[(I-1)] The innovations $(\xi_{t})_{t\in\Z}$
satisfy $m_p:=\E[\llvert \xi\rrvert ^p ]<\infty$.

\item[(I-2)] The innovations $(\xi_{t})_{t\in
\Z}$ satisfy
$\phi(\zeta):=\E[\rme^{\zeta\llvert \xi\rrvert } ]<\infty$.
\end{longlist}
The following one will be
used to obtain a lower bound.
\begin{longlist}[(I-3)]\label{hyp:innov-lowerbound}
\item[(I-3)] The innovations $(\xi_{t})_{t\in\Z}$ admit a density $f$ such that
\[
\kappa= \sup_{v\neq0}v^{-2}\int f (u )\log{
 \frac{f (u )}{f (u+v )}} \,\rmd u < \infty . %
\]
\end{longlist}

Assumption~(I-3) is standard for proving lower
bounds in
nonparametric regression estimation, see
\citet{Tsybakov2009}, Chapter~2. It is
satisfied by Gaussian density with $\kappa=1$.

\subsubsection{Nonparametric setting}
The setting of Definition~\ref{definition:tvar} and of assumptions derived
thereafter is essentially nonparametric, since for given initial distribution~$\psi$, the distribution of the observations $X_{1,T},\dots,X_{T,T}$ are
determined by the unknown parameter function $(\btheta,\sigma)$. The doubly
indexed $X_{t,T}$ refers to the fact that this distribution cannot be
seen as a
distribution on $\R^\Z$ marginalized on $\R^T$ as the usual time
series setting
but rather as a sequence of distributions on $\R^T$ indexed by
$T$. It corresponds to the usual nonparametric approach for
studying statistical inference based on this model. In this
contribution, we
focus on the prediction problem, which is to answer the question: for given
smoothness conditions on $(\btheta,\sigma)$, what is the mean
prediction error
for predicting $X_{t,T}$ from its past? The standard nonparametric approach
is to answer this question in a minimax sense by determining, for a given
sequence of predictors $\widehat{X}_T=(\widehat{X}_{t,T})_{1\leq
t\leq T}$, the
maximal risk
%
\begin{eqnarray}\label{eq:max-risk}
&& S_T (\widehat{X}_T;\psi,\beta,R,\delta,
\rho,\sigma_+ )
\nonumber\\[-8pt]\\[-8pt]\nonumber
&&\qquad = \sup_{(\btheta,\sigma)} { \frac{1}{T}}
\sum_{t=1}^{T} \biggl(\E^{\psi
}_{ (\btheta,\sigma)}
\bigl[ (\widehat{X}_{t,T}-X_{t,T} )^{2} \bigr] -
\sigma^{2} \biggl({ \frac{t}{T}} \biggr) \biggr) ,
\end{eqnarray}
where:
\begin{longlist}[(a)]
\item[(a)] $\widehat{X}_T$ is assumed to belong to $\mathcal{P}_T$ as in
Definition~\ref{def:predictor},

\item[(b)] the sup is taken over $(\btheta,\sigma)\in\mathcal{C} (\beta
,R,\delta,\rho,\sigma_{+} )$ within a smoothness class of
functions,

\item[(c)] the expectation $\E^{\psi}_{ (\btheta,\sigma)}$
is that
associated to assumption~(M-2).
\end{longlist}
The reason for subtracting the average $\sigma^2(t/T)$ over all
$1\leq t\leq T$ in this prediction risk is that it corresponds to the best
prediction risk, would the parameters $(\btheta,\sigma)$ be exactly known.
We observe that dividing $X_{t,T}$ by the class parameter $\sigma_+$
amounts to
take $\sigma_+=1$. In addition, we have
\[
S_T (\widehat{X}_T;\psi,\beta,R,\delta,\rho,\sigma_+ )=
\sigma_{+}^2 S_T \bigl(\widehat{X}_T
\sigma_{+}^{-1};\psi,\beta,R,\delta,\rho,1 \bigr) ,
\]
so the prediction problem in the class
$\mathcal{C} (\beta,R,\delta,\rho,\sigma_{+} )$ can be
reduced to the prediction problem in the class
$\mathcal{C} (\beta,R,\delta,\rho,1 )$. Accordingly, we
define the
reduced minimax risk by
%
\begin{eqnarray}
\label{eq:reduced-minimax-risk}
&& \overline{M}_T (\psi,\beta,R,\delta ,\rho)\nonumber
\\
&&\qquad = \inf _{\widehat{X}_T\in\mathcal{P}_T}S_T (\widehat{X}_T;\psi,\beta,R,
\delta,\rho,1 )
\\
&&\qquad =\inf_{\widehat{X}_T\in\mathcal{P}_T}\sigma_+^{-2}
 S_T (\widehat{X}_T;\psi,\beta,R,\delta,\rho,\sigma_+ )
\qquad\mbox{for all $\sigma_+>0$}  .\nonumber
\end{eqnarray}

In Section~\ref{section:lower_bound}, we provide a lower bound of the
minimax rate
in the case where the smoothness class is of the form
$\mathcal{C} (\beta,R,\delta,\rho,\sigma_{+} )$. Then, in
Section~\ref{subsection:minimax}, relying on the aggregation oracle
bounds of
Section~\ref{subsection:general:bound}, we derive an upper bound with
the same
rate as the lower bound using the same smoothness class of the
parameters. Moreover, we exhibit an online predictor which does not
require any knowledge
about the smoothness class and which is thus minimax adaptive. In other
words, it
is able to adapt to the unknown smoothness of the parameters from the
data. To
our knowledge, such theoretical results are new for locally stationary models.

\subsection{Lower bound} \label{section:lower_bound}

A lower bound on the minimax rate for the estimation error of $\btheta
$ is given by
\citeauthor{MoulinesPriouretRoueff2005} [(\citeyear
{MoulinesPriouretRoueff2005}), Theorem~4]. Clearly, a~predictor
\[
\widehat X_{t,T}=\sum_{k=1}^d
\widehat{\btheta}_{t,T}(k) X_{t-k,T}
\]
can be defined from an estimator $\widehat{\btheta}_{t,T}$, and the resulting
prediction rate can be controlled using the estimation rate (see
\citeauthor{GiraudRoueffSanchez-Perez2015s} [(\citeyear
{GiraudRoueffSanchez-Perez2015s}), Section~B.1] for the details). The
next theorem
provides a lower bound of
the minimax rate of the risk of \emph{any} predictor of the process
$(X_{t,T})_{1\leq t\leq T}$.
Combining this result with [\citeauthor{GiraudRoueffSanchez-Perez2015s}
(\citeyear{GiraudRoueffSanchez-Perez2015s}), Lemma~9], we
show that a predictor
obtained by [Giraud, Roueff and Sanchez-Perez (\citeyear
{GiraudRoueffSanchez-Perez2015s}), equation~(B.1)] from a
minimax rate estimator of
$\btheta$ automatically achieves the minimax
prediction rate.

%
\begin{theorem}\label{theorem:lower_bound_quadratic_risk}
Let $\delta\in(0,1)$, $\beta>0$, $R>0$ and $\rho\in[0,1]$. Suppose that
assumption~\textup{(M-2)} holds and assume \textup{(I-3)}
on the distribution $\psi$ of the
innovations. Then we have
%
\begin{equation}
\label{eq:lower-bound-minimax} \mathop{\operatorname{lim}\operatorname{inf}}_{T\to\infty
}T^{{2\beta}/(1+2\beta)}
\overline{M}_T(\psi,\beta,R,\delta,\rho)>0  ,
\end{equation}
where $\overline{M}_T$ is defined in~(\ref{eq:reduced-minimax-risk}).
\end{theorem}

The proof is postponed to
Section~\ref{section:proof_lower_bound}.

\subsection{Minimax adaptive forecasting of the TVAR
process} \label{subsection:minimax}

In \citet{Arkoun2011}, an adaptive estimator of the autoregressive
function of
a Gaussian TVAR process of order 1 is studied. It relies on the
Lepski{\u\i}'s procedure [see \citet{Lepskii1990}], which seems
difficult to
implement in an online context.

Our minimax adaptive predictor is based on the aggregation of
sufficiently many
predictors, assuming that at least one of them converges at the minimax
rate. The oracle bounds found in Section~\ref
{subsection:general:bound} imply
that the aggregated predictor is minimax rate adaptive under appropriate
assumptions. Seminal works using the
aggregation to adapt to the minimax convergence rate are \citet
{Yang2000} (nonparametric regression) and
\citet{Yang-TAS2000} (density estimation); see also \citet
{Catoni2004} for a more general presentation.

In the TVAR model~(M-2), it is natural to consider $L$-Lipschitz
predictors $(\widehat{X}_{t,T})_{1\leq t\leq T}$ of
$(X_{t,T})_{1\leq t\leq T}$ with a sequence $L$ supported on
$\{1,\dots,d\}$. Then $L^*$ in~(\ref{eq:hyp:lip}) corresponds to the maximal
$\ell^1$-norm of the TVAR parameters. Since for the process itself to be
stable, this norm has to be bounded independently of~$T$,
condition~\textup{(L-1)} is a
quite natural assumption for the TVAR model; see
\citeauthor{GiraudRoueffSanchez-Perez2015s} [(\citeyear
{GiraudRoueffSanchez-Perez2015s}), Section~B.1] for the details.

A practical advantage of the proposed procedures is that, given a set of
predictors that behaves well under specific smoothness assumptions, we obtain
an aggregated predictor which performs almost as well as or better than the
best of these predictors, hence which behaves well without any prior
knowledge on
the smoothness of the unknown parameter. Such an adaptive property can be
formally demonstrated by exhibiting an adaptive minimax rate for the aggregated
predictor which coincides with the lower bound given in
Theorem~\ref{theorem:lower_bound_quadratic_risk}.

The first ingredient that we need is the following.

%
\begin{definition}[{[$(\psi,\beta)$-minimax-rate predictor]}]\label{def:minimax-rate-expert}
Let $\psi$ be a distribution on~$\R$ and $\beta>0$. We say that
$\widehat{X}=(\widehat{X}_{T})_{T\geq1}$ is a $(\psi,\beta
)$-minimax-rate sequence of predictors if, for all $T\geq1$,
$\widehat{X}_T\in\mathcal{P}_T$ and, for all $\delta\in(0,1)$, $R>0$,
$\rho\in(0,1]$ and $\sigma_+>0$,
%
\begin{equation}
\label{equation:definition_minimax} \limsup_{T\to\infty}T^{{2\beta
}/(1+2\beta)} S_T
(\widehat{X}_{T};\psi,\beta,R,\delta,\rho,\sigma_+ )<\infty  ,
\end{equation}
where $S_T$ is defined by~(\ref{eq:max-risk}).
\end{definition}

The term \emph{minimax-rate} in this definition refers to the fact
that the
maximal rate in~(\ref{equation:definition_minimax}) is equal to the
minimax lower bound~(\ref{eq:lower-bound-minimax}) for the class
$\mathcal{C} (\beta,R,\delta,\rho,\sigma_{+} )$.
We explain in \citeauthor{GiraudRoueffSanchez-Perez2015s} [(\citeyear
{GiraudRoueffSanchez-Perez2015s}), Section~B] how to build such
predictors which
are moreover
$L$-Lipschitz for some $L$ only depending on $d$.
To adapt to an unknown smoothness, we rely on a collection of
$(\psi,\beta)$-minimax-rate predictors with $\beta$ within $(0,\beta
_0)$, where
$\beta_0$ is the (possibly infinite) maximal smoothness index.

%
\begin{definition}[(Locally bounded set of $\psi$-minimax-rate predictors)]
\label{def:minimax-rate-expert-loc-bounded}
Let $\psi$ be a distribution on $\R$ and $\beta_{0}\in(0,\infty]$.
We say that
$\{\widehat{X}^{(\beta)},  \beta\in(0,\beta_0)\}$ is a locally
bounded set of $\psi$-minimax-rate predictors if for each $\beta$,
$\widehat{X}^{(\beta)}$ is a $(\psi,\beta)$-minimax-rate predictor
and if
moreover, for all $\delta\in(0,1)$, $R>0$, $\rho\in(0,1]$, $\sigma
_+>0$ and for each closed interval
$J\subset(0,\beta_0)$,
\[
\limsup_{T\to\infty} \sup_{\beta\in J}
T^{{2\beta}/(1+2\beta)} S_T \bigl(\widehat{X}_{T}^{(\beta)};
\psi,\beta,R,\delta,\rho,\sigma_+ \bigr) <\infty , %
\]
where $S_T$ is defined by~(\ref{eq:max-risk}).
\end{definition}

The following lemma shows that, given a locally bounded set of minimax-rate
predictors, we can always pick a finite subset of at most $N=\lceil
(\log
T)^2\rceil$ predictors among which the best one achieves the minimax
rate of
any unknown smoothness index.

%
\begin{lemma}\label{lemma:minimax_experts}
Let $\psi$ be a distribution on $\R$.
Let $\beta_0\in(0,\infty]$ and $\{\widehat{X}^{(\beta)},
\beta\in(0,\beta_0) \}$ be a corresponding locally bounded set of
$\psi$-minimax-rate predictors. Set, for any $N\geq1$,
%
\begin{equation}
\label{eq:betai-def} \beta_i= \cases{ (i-1)\beta_{0}/N, &\quad
if $\beta_0<\infty$,
\vspace*{3pt}\cr
(i-1)/N^{1/2}, &\quad otherwise,}
\quad1\leq i\leq N  .
\end{equation}
Suppose moreover, in the case where $\beta_0<\infty$, that $N\geq
\lceil\log T\rceil$, and, in the case where $\beta_0=\infty$, that
$N\geq
\lceil(\log T )^{2} \rceil$.
Then we have, for all $\beta\in(0,\beta_{0} )$,
$\delta\in(0,1)$, $R>0$, $\rho> 0$ and $\sigma_+>0$,
\[
\limsup_{T\to\infty}T^{{2\beta}/(1+2\beta)} \min_{i=1,\dots,N}S_T
\bigl(\widehat{X}_{T}^{(\beta_i)};\psi,\beta,R,\delta,\rho,
\sigma_+ \bigr)<\infty . %
\]
\end{lemma}

The proof of this lemma is postponed to
\citeauthor{GiraudRoueffSanchez-Perez2015s} [(\citeyear
{GiraudRoueffSanchez-Perez2015s}), Section~A.8].
Lemma~\ref{lemma:minimax_experts} says that to obtain a minimax-rate
predictor which adapts to an unknown smoothness index $\beta$, it is
sufficient to select it judiciously among $\log T$ or $(\log T)^2$ well chosen
nonadaptive minimax-rate predictors.

As a consequence of Theorem~\ref{thm:oracle-bounds} and Lemma~\ref
{lemma:minimax_experts}, we obtain an adaptive predictor by aggregating
them (instead of selecting one of them), as stated in the following
result.

%
\begin{theorem}\label{theorem:upper_bound_aggregated_forecaster_adaptative}
Let $\psi$ be a distribution on $\R$.
Let $\beta_0\in(0,\infty]$ and
$\{\widehat{X}^{(\beta)}$, $ \beta\in(0,\beta_0)\}$ be a locally
bounded\vspace*{1pt} set of $\psi$-minimax-rate and $L$-Lipschitz predictors with $L$
satisfying~\textup{(L-1)}. Define $(\widehat{X}_{t,T})_{1\leq t\leq
T}$ as
the predictor aggregated from $\{\widehat{X}^{(\beta_i)}, 1\leq
i\leq N
\}$ with $N$ defined by
%
\begin{equation}
\label{eq:def-N-agreg-tvar} N= \cases{ \lceil\log T\rceil, &\quad if
$\beta_0<
\infty$,
\vspace*{3pt}\cr
\bigl\lceil(\log T)^2\bigr\rceil, &\quad otherwise,}
\end{equation}
$\beta_i$ defined by~(\ref{eq:betai-def}), and with weights defined
according to one of the following setting depending on the
assumption on $\psi$ and $\beta_0$:
\begin{longlist}
\item[(i)]  If\vspace*{1pt} $\psi$ satisfies~\textup{(I-1)} with
$p\geq4$ and $\beta_0\leq1/2$, use the
weights~(\ref{definition:alphas_gradient_quadratic_loss}) with $\eta=
\sigma_+^{-2}(\log(\lceil\log T\rceil)/T )^{1/2}$.

\item[(ii)]  If $\psi$ satisfies~\textup{(I-1)} with
$p>2$ and $\beta_0\leq(p-2)/4$, use the
weights~(\ref{definition:alphas_quadratic_loss}) with $\eta=
\sigma_+^{-2}(\log(\lceil\log T\rceil)/T)^{2/p}$.

\item[(iii)]  If $\psi$ satisfies~\textup{(I-2)}, use the
weights~(\ref{definition:alphas_quadratic_loss})
with $\eta= \sigma_{+}^{-2} (\log
T )^{-3}$.
\end{longlist}
Then we have, for any $\beta\in(0,\beta_0 )$,
$\delta\in(0,1)$, $R>0$, $\rho\in(0,1]$ and $\sigma_+>0$,
%
\begin{equation}
\label{eq:final-tvar-strat12} \limsup_{T\to\infty}T^{{2\beta}/(1+2\beta
)} S_T
(\widehat{X}_{T};\psi,\beta,R,\delta,\rho,\sigma_+ )<\infty .
\end{equation}
\end{theorem}

The proof of this theorem is postponed to Giraud, Roueff and Sanchez-Perez 
[(\citeyear{GiraudRoueffSanchez-Perez2015s}), Section~A.9].

%
\begin{remark} \label
{remark:theorem_upper_bound_aggregated_forecaster_adaptative}
The\vspace*{1pt} limitation to $\beta_0\leq1/2$ in~(i) under
assumption~(I-1) for $\psi$ follows from the factor
$(\log N/T)^{1/2}$ obtained in the oracle inequality~(\ref{eq:strat1-final})
of Theorem~\ref{thm:oracle-bounds} after optimizing in $\eta$
[see~(\ref{eq:strat1-final-eta-opt})]. If $p>4$ this
restriction is weakened to $\beta_0\leq(p-2)/4$ in~(ii)
taking into account the factor
$(\log N/T)^{1-2/p}$ obtained in the oracle inequality~(\ref
{eq:strat2_p-final})
of Theorem~\ref{thm:oracle-bounds} after optimizing in $\eta$
[see~(\ref{eq:strat2_p-final-eta-opt})]. In the last case, the
limitation of $\beta_{0}$ drops when applying the oracle
inequality~(\ref{eq:strat2_exp-final})
of the same theorem. However, a stronger condition on $\psi$ is then
required.
\end{remark}

%
\begin{remark}\label{rem:beta0-eta-minim-adapt-forec}
It may happen that the locally bounded set of $\psi$-minimax-rate predictors
is limited to some $\beta_{0}<\infty$ [see the example of the NLMS predictors
in \citet{GiraudRoueffSanchez-Perez2015s}, Section~B.2]. In this case, the result
roughly needs $\log T$ predictors and the computation of the aggregated one
requires less operations than if $\beta_{0}$ were infinite. For these
reasons, we do not consider in general that $\beta_{0}=\infty$. On
the one
hand, a finite
$\beta_0$ yields a restriction on the set of (unknown) smoothness
indices $\beta$ for which the aggregated predictors are minimax rate
adaptive. On the other hand, if $\beta_0=\infty$,
Theorem~\ref{theorem:upper_bound_aggregated_forecaster_adaptative} then
requires the stronger assumption~(I-2) on the
process.
\end{remark}

%
\begin{remark}\label{rem:sigma-eta-minim-adapt-forec}
The constant $\sigma_{+}^{-2}$ present in the definitions of $\eta$
in the
three cases~(i), (ii) and (iii) corresponds
to the
homogenization of the remaining terms appearing in
Theorem~\ref{thm:oracle-bounds} [the second lines of~(\ref{eq:strat1-final}),
(\ref{eq:strat2_p-final}) and (\ref{eq:strat2_exp-final})]. Indeed with
the proposed choices and in the three cases, the constant $\sigma_{+}^{2}$
factors out in front of the remaining terms [see the last three displayed
equations in \citet{GiraudRoueffSanchez-Perez2015s}, Section~A.9]. However,
the $\sigma_+^{-2}$
in the definitions of $\eta$ does not impact the convergence rate in
the sense
that Theorem~\ref
{theorem:upper_bound_aggregated_forecaster_adaptative} is
still valid using any other constant ($1$, e.g.) in these definitions.
\end{remark}

\section{Proofs of the upper bounds}
\label{subsection:oracle_inequalities}
\subsection{Preliminary results}
We start with a lemma which gathers useful adaptations of well-known
inequalities
applying to the aggregation of deterministic predicting sequences.

%
\begin{lemma}\label{lemma:deterministic_upper_bound_aggregated_forecaster}
Let $(x_{t})_{1\leq t\leq T}$ be a real valued sequence and
$\{(\widehat{x}_{t}^{(i)})_{1\leq t\leq T}, 1\leq i \leq N\}$ be a
collection of predicting sequences. Define $(\widehat{x}_{t})_{1\leq
t\leq T}$ as the sequence of aggregated predictors obtained from this
collection with the
weights (\ref{definition:alphas_gradient_quadratic_loss}). Then, for any
$\eta>0$, we have
%
%
\begin{eqnarray}
\label{eq:strat1-deterministic}
\frac{1}{T}\sum_{t=1}^T
(\widehat{x}_{t}-x_{t} )^2 &\leq& \inf
_{\bolds{\nu}\in
\mathcal{S}_N}\frac{1}{T}\sum_{t=1}^T
\bigl(\widehat{x}_{t}^{[\bolds
\nu]}-x_{t}
\bigr)^2 + \frac{\log N}{T\eta} + \frac{2\eta}{T}\sum
_{t=1}^{T} y_{t}^{4} ,
\end{eqnarray}
where $y_{t}=\llvert x_{t}\rrvert +\max_{1\leq i\leq
N}\llvert \widehat{x}_{t}^{(i)}\rrvert $.

Define now $(\widehat{x}_{t})_{1\leq
t\leq T}$ as the sequence of aggregated predictors obtained with the
weights~(\ref{definition:alphas_quadratic_loss}). Then, for any $\eta
>0$, we have
\begin{eqnarray}\label{eq:strat2-deterministic-allcases}
&& \frac{1}{T}\sum_{t=1}^T
(\widehat{x}_{t}-x_{t} )^2
\nonumber\\[-8pt]\\[-8pt]\nonumber
&&\qquad \leq\min
_{i=1,\dots,N}\frac{1}{T}\sum_{t=1}^T
\bigl(\widehat{x}_{t}^{(i)}-x_{t}
\bigr)^2 + \frac{\log N}{T\eta} +\frac{1}{T}\sum
_{t=1}^{T} \biggl(y_{t}^{2}-
\frac{1}{2\eta} \biggr)_{+} ,
\end{eqnarray}
where $y_{t}=\llvert x_{t}\rrvert +\max_{1\leq i\leq
N}\llvert \widehat{x}_{t}^{(i)}\rrvert $.
\end{lemma}

\begin{pf}
With the weights defined by~(\ref
{definition:alphas_gradient_quadratic_loss}), by
slightly adapting [\citeauthor{Stoltz2011} (\citeyear{Stoltz2011}),
Theorem 1.7], we have that
\begin{eqnarray*}
\frac{1}{T}\sum_{t=1}^T (
\widehat{x}_{t}-x_{t} )^2-\inf
_{\bolds{\nu}\in\mathcal{S}_N}\frac{1}{T}\sum_{t=1}^T
\bigl(\widehat{x}_{t}^{[\bolds\nu]}-x_{t}
\bigr)^2 &\leq& \frac{\log
N}{T\eta} + \frac{\eta}{8T}
 s_{T}^{*} ,
\end{eqnarray*}
where $s_{T}^{*}=\sum_{t=1}^{T}s_{t}^{2}$ and $s_{t}=2\max_{1\leq
i\leq N}\llvert 2(\sum_{j=1}^{N}\widehat{\alpha}_{j,t}\widehat
{x}_{t}^{(j)}-x_{t})\widehat{x}_{t}^{(i)}\rrvert $.
The bound (\ref{eq:strat1-deterministic}) follows by using that
$ \{\widehat{\alpha}_{i,t} \}_{1\leq i\leq N}$ is in the
simplex $\mathcal{S}_N$
defined in~(\ref{eq:def-simplexN}).

We now
prove~(\ref{eq:strat2-deterministic-allcases}). We adapt the proof of
\citeauthor{Catoni2004} [(\citeyear{Catoni2004}), Proposition 2.2.1.]
to unbounded sequences by
replacing the
convexity argument by the following lemma.

\begin{lemma} \label{lemma:bound_integral_concave_exponential}
Let $a>0$ and $\PP$ a probability distribution supported on $[-a,a]$.
Then we have
\begin{eqnarray*}
&&\int\exp\bigl(-x^{2} \bigr)\,\rmd\PP(x ) \leq\exp\biggl(- \biggl(
\int x\,\rmd\PP(x ) \biggr)^{2}+ \biggl(a^{2}-
\frac{1}{2} \biggr)_{+} \biggr)  .
\end{eqnarray*}
\end{lemma}

The proof of Lemma~\ref{lemma:bound_integral_concave_exponential} is
postponed to
Section~\ref{sec:proof-lemma-technical-exp-concave} in \hyperref[section:useful_bounds_TVAR_postponed_proofs]{Appendix}.
Now, let $\eta>0$ and $t=1,\ldots,T$. Using
Lemma~\ref{lemma:bound_integral_concave_exponential} with the probability
distribution $\PP$ defined by
$\PP(A)=\sum_{i=1}^{N}\widehat{\alpha}_{i,t}\1_A(\eta
^{1/2}(\widehat{x}_{t}^{(i)}-x_{t}))$
and $a=\eta^{1/2}y_{t}$,
we get that
\[
\sum_{i=1}^{N}\widehat{
\alpha}_{i,t}\exp\bigl(-\eta\bigl(\widehat{x}_{t}^{(i)}-x_{t}
\bigr)^{2} \bigr) \leq\exp\biggl(-\eta(\widehat{x}_{t}-x_{t}
)^{2}+\eta\biggl(y_{t}^{2}-\frac{1}{2\eta}
\biggr)_{+} \biggr) . %
\]
Taking the log, multiplying by $-\eta^{-1}$ and re-ordering the terms,
we obtain that
\begin{eqnarray*}
(\widehat{x}_{t}-x_{t} )^{2} &\leq&-
\frac{1}{\eta}\log\Biggl(\sum_{j=1}^{N}
\widehat{\alpha}_{i,t}\exp\bigl(-\eta\bigl(\widehat{x}_{t}^{(i)}-x_{t}
\bigr)^{2} \bigr) \Biggr)
\\
&&{} + \biggl(y_{t}^{2}-
\frac{1}{2\eta} \biggr)_{+}  . %
\end{eqnarray*}
Taking the average over $t=1,\dots,T$ and developing the expression of
$\widehat{\alpha}_{i,t}$, we obtain
\begin{eqnarray}\label{equation:regret_deterministic_str2}
\frac{1}{T}\sum_{t=1}^{T}
(x_{t}-\widehat{x}_{t} )^{2} &\leq&-
\frac{1}{\eta T}\log\Biggl(\frac{1}{N}\sum_{i=1}^{N}
\exp\Biggl(-\eta\sum_{t=1}^{T} \bigl(
\widehat{x}_{t}^{(i)}-x_{t} \bigr)^{2}
\Biggr) \Biggr)
\nonumber\\[-8pt]\\[-8pt]\nonumber
&&{} +\frac{1}{T}\sum_{t=1}^{T}
\biggl(y_{t}^{2}-\frac{1}{2\eta} \biggr)_{+}.
\end{eqnarray}
Using that
$\sum_{i=1}^{N}\exp(-\eta\sum_{t=1}^{T}(\widehat
{x}_{t}^{(i)}-x_{t})^{2})\geq
\exp(-\eta\min_{i=1,\ldots,N}\sum_{t=1}^{T}(\widehat
{x}_{t}^{(i)}-x_{t})^{2})$, we get the bound~(\ref
{eq:strat2-deterministic-allcases}).
\end{pf}

\subsection{Proof of Theorem~\protect\ref{thm:oracle-bounds}}
\label{sec:proof-theor-strat12}

We prove the cases~(i), (ii) and~(iii) successively. We denote $Y_{t}=\llvert
{X}_{t}\rrvert +\max_{1\leq
i\leq N}\llvert \widehat{X}_{t}^{(i)}\rrvert $.

\begin{longlist}[Case (iii)]
\item[\textit{Case}~(i).]
Applying~(\ref{eq:strat1-deterministic}) in
Lemma~\ref{lemma:deterministic_upper_bound_aggregated_forecaster}
with $\E[\inf\cdots]\leq\inf\E[\cdots]$, we obtain
%
\begin{eqnarray}
\label{inequality:expectation_aggregated_forecaster_moment} \frac
{1}{T}\sum_{t=1}^T
\E\bigl[ (\widehat{X}_{t}-X_{t} )^2 \bigr]
&\leq&\inf_{\bolds{\nu}\in\mathcal{S}_N}\frac{1}{T}\sum
_{t=1}^T\E\bigl[ \bigl(\widehat{X}_{t}^{[\bolds\nu]}-X_{t}
\bigr)^2 \bigr]
\nonumber\\[-8pt]\\[-8pt]\nonumber
&&{} + \frac{\log N}{T\eta} + \frac{2\eta}{T}\sum
_{t=1}^{T} \E\bigl[Y_{t}^{4} \bigr]  .
\end{eqnarray}
Using that the predictors are $L$-Lipschitz and the process
$(X_t)_{t\in\Z}$ satisfies~(M-1), we have, for all $1\leq
t\leq T$,
\begin{eqnarray}
\label{inequality:X-Z}
Y_{t} &=& \llvert{X}_{t}\rrvert+\max
_{1\leq i\leq N}\bigl\llvert\widehat{X}_{t}^{(i)}
\bigr\rrvert\nonumber
\\
&\leq& \sum_{j\in\Z}A_{t}(j)  Z_{t-j}+\sum_{s\geq1}\sum
_{j\in\Z}L_s A_{t-s}(j)
 Z_{t-s-j}
\\
&\leq& \sum_{j\in\Z} B_t(j) Z_{t-j},\nonumber
\end{eqnarray}
where
\[
B_t(j)=A_{t}(j)+\sum_{s\geq1}L_{s}
 A_{t-s}(j-s) .
\]
Applying the Minkowski inequality together with~(\ref
{inequality:X-Z}), (\ref{eq:unif-coeff-bound}) and~(\ref
{eq:hyp:lip}), we obtain, for all $1\leq t\leq T$,
\begin{eqnarray*}
&&\E\bigl[Y_{t}^{4} \bigr] \leq\E\biggl[ \biggl(\sum
_{j\in\Z} B_t(j) Z_{t-j}
\biggr)^{4} \biggr] \leq A_{*}^4(1+L_*)^{4}
 \sup_{t\in
\Z}\E\bigl[Z_t^4 \bigr] .
\end{eqnarray*}
Since the process $Z$ fulfills~\textup{(N-1)} with $p=4$,
plugging this bound in~(\ref
{inequality:expectation_aggregated_forecaster_moment}) we obtain~(\ref
{eq:strat1-final}).

\item[\textit{Case}~(ii).]
We use~(\ref
{eq:strat2-deterministic-allcases}) in
Lemma~\ref{lemma:deterministic_upper_bound_aggregated_forecaster} and the
inequality $(x^2-1/(2\eta))_{+}\leq(2\eta)^{p/2-1}x^{p}$ which holds for
$x\geq0$ and $p\geq2$. We get, taking the expectation,
\begin{eqnarray}
\label{inequality:expectation_aggregated_forecaster_moment_1}
\qquad\frac{1}{T}\sum_{t=1}^{T}
\E\bigl[ (\widehat{X}_{t,T}-X_{t,T} )^2 \bigr]
&\leq&\min_{i=1,\dots,N}\frac{1}{T}\sum
_{t=1}^{T}\E\bigl[ \bigl(\widehat{X}_{t,T}^{(i)}-X_{t,T}
\bigr)^2 \bigr] + \frac{\log
N}{T\eta}
\nonumber\\[-8pt]\\[-8pt]\nonumber
&&{} +(2\eta)^{p/2-1}\max_{t=1,\ldots,T}\E\bigl[Y_{t}^p
\bigr]  .
\end{eqnarray}
Applying the Minkowski inequality, (\ref{inequality:X-Z}) and
assumption~\textup{(N-2)},
\[
\E\bigl[Y_{t}^{p} \bigr] \leq\biggl(\sum
_{j\in Z}B_{t}(j) \bigl(\E\bigl[Z_{t-j}^{p}
\bigr] \bigr)^{1/p} \biggr)^{p} \leq A_{*}^{p}(1+L_*)^{p}
 \sup_{t\in Z}\E\bigl[Z_t^{p} \bigr] .
\]
Using this bound which is independent of $t$, with~\textup{(N-1)}
and~(\ref{inequality:expectation_aggregated_forecaster_moment_1}), the
inequality~(\ref{eq:strat2_p-final}) follows.

\item[\textit{Case}~(iii).]
To obtain~(\ref{eq:strat2_exp-final}), we
again use~(\ref{eq:strat2-deterministic-allcases}) in
Lemma~\ref{lemma:deterministic_upper_bound_aggregated_forecaster} but
now with
an exponential bound for $(Y_t^{2}-1/(2\eta))_{+}$. We note that, or
all $u>0$,
\[
\sup_{x\geq1} \bigl(x^2-1 \bigr)\rme^{-ux}=
\bigl(x_0^2-1 \bigr)\rme^{-u x_0}\qquad\mbox{with
} x_0=u^{-1} \bigl(1+ \bigl(1+u^2
\bigr)^{1/2} \bigr)  . %
\]
It follows that, for all $x\in\R$ and $u>0$,
\[
\bigl(x^2-1 \bigr)_+\leq\rme^{ux} \bigl(x_0^2-1
\bigr)\rme^{-u
x_0}\leq\rme^{ux}2u^{-2} (2+u )
\rme^{-1- u}  . %
\]
Applying this bound with $x=(2\eta)^{1/2}Y_t$ and $u=\lambda(2\eta
)^{-1/2}$ we
get
\[
\biggl(Y_t^2-\frac{1}{2\eta} \biggr)_{+}=(2
\eta)^{-1} \bigl(x^2-1 \bigr)_{+} \leq2
\lambda^{-2} \bigl(2+\lambda(2\eta)^{-1/2} \bigr)
\rme^{-1-
\lambda(2\eta)^{-1/2}} \rme^{\lambda Y_t} . %
\]
Plugging this into~(\ref{eq:strat2-deterministic-allcases}) and taking the
expectation, we obtain that
\begin{eqnarray}\label{inequality:expectation_aggregated_forecaster_exponential_moment}
&& \frac{1}{T}\sum_{t=1}^{T}
\E\bigl[ (\widehat{X}_{t,T}-X_{t,T} )^2 \bigr]\nonumber
\\
&&\qquad \leq\min_{i=1,\dots,N}\frac{1}{T}\sum
_{t=1}^{T}\E\bigl[ \bigl(\widehat{X}_{t,T}^{(i)}-X_{t,T}
\bigr)^2 \bigr] + \frac{\log
N}{T\eta}
\\
&&\quad\qquad{} +2\lambda^{-2} \bigl(2+\lambda(2\eta)^{-1/2} \bigr)
\rme^{-1- \lambda(2\eta)^{-1/2}}  \max_{t=1,\ldots,T}\E\bigl[\rme^{\lambda Y_t}
\bigr]  .\nonumber
\end{eqnarray}
We now use assumption~\textup{(N-2)}.
Since $B_t(j)\leq a^*(1+L_{*})$ for all $j,t\in\Z$ and
\[
\sum_{j\in\Z}B_t(j)\leq
A_{*}(1+L_{*}),
\]
Jensen's inequality and (\ref{inequality:X-Z}) gives that, for any
$\lambda\leq\zeta/(a^*(1+L_{*}))$,
\begin{eqnarray*}
\E\bigl[\rme^{\lambda Y_{t}} \bigr]&\leq&\E\bigl[\rme^{\lambda
(\llvert {X}_{t}\rrvert +\max_{1\leq i\leq N}\llvert \widehat
{X}_{t}^{(i)}\rrvert ) } \bigr]
\\
&\leq&\prod_{j\in\Z}\E\bigl[\rme^{\lambda B_t(j) Z_{t-j}} \bigr]
\\
&\leq&\prod_{j\in\Z} \bigl(\phi(\zeta)
\bigr)^{\lambda
B_t(j)/\zeta} \leq\bigl(\phi(\zeta) \bigr)^{\lambda
A_{*}(1+L_{*})/\zeta}.
\end{eqnarray*}

The combination of this bound with~(\ref
{inequality:expectation_aggregated_forecaster_exponential_moment})
gives~(\ref{eq:strat2_exp-final}). The proof of Theorem~\ref
{thm:oracle-bounds} is complete.
\end{longlist}

\subsection{Proof of case~\textup{(iii)} in Corollary~\protect
\ref{corollary:oracle-bounds}}
\label{sec:proof-corollary-strat12}

Minimizing the sum of the two terms appearing in the second line
of~(\ref{eq:strat2_exp-final}) is a bit more involved, since it
depends both on
$\eta$ and $\lambda$. Under condition~(\ref{eq:def-sup-unif-coeff-bound}), the
quantity $ (\mphi(\zeta) )^{\lambda A_*(1+L_*)/\zeta}$ remains
between two positive constants while, for any $\eta>0$,
$\lambda^{-2}(2+\lambda(2\eta)^{-1/2})$ is decreasing as $\lambda$
increases. To simplify $ (\mphi(\zeta) )^{\lambda
A_*(1+L_*)/\zeta}$
into $\mphi(\zeta)$, we simply take
\[
\lambda=\frac{\zeta}{A_{*}(1+L_*)} ,
\]
which
satisfies~(\ref{eq:def-sup-unif-coeff-bound}). Now that $\lambda$ is
set, it remains to choose a
value of $\eta$ which (almost) minimizes
\[
\frac{\log N}{T\eta} + \frac{2\mphi(\zeta)}{\rme}\lambda^{-2} \bigl
(2+\lambda(2
\eta)^{-1/2} \bigr) \rme^{- \lambda(2\eta)^{-1/2}} . %
\]
The $\eta$ defined as
in~(\ref{eq:strat2_exp-choice-eta}) is chosen so that $(\log N)/T=\rme^{-
\lambda(2\eta)^{-1/2}}$, and we get~(\ref{eq:strat2_exp-final-opt-eta}).

\section{Proof of the lower bound} \label{section:proof_lower_bound}
We now provide a proof of
Theorem~\ref{theorem:lower_bound_quadratic_risk}. We consider an
autoregressive equation of order one
%
\begin{equation}
\label{eq:lb-tvar-eq} X_{t,T}=\theta\biggl(\frac{t-1}{T} \biggr)
X_{t-1,T}+\xi_{t},
\end{equation}
where $(\xi_{t})_{t\in\Z}$ is i.i.d. with density $f$ as
in~(I-3). In this case, if $\sup_{u\leq
1}\llvert \theta(u)\rrvert <1$, the
representation~(\ref{eq:inf-series-representation}) of the stationary
solution reads, for all $ t\leq T$ as
%
\begin{equation}
\label{eq:inf-series-representation-d1} X_{t,T} = \sum_{j=0}^{\infty}
\prod_{s=1}^{j}\theta\biggl(
\frac
{t-s}{T} \biggr)  \xi_{t-j} ,
\end{equation}
with the convention $ \prod_{s=1}^{0}\theta((t-s)/T)=1$.
The class of models so defined with $\theta\in\Lambda_{1}(\beta,R)\cap
s_{1}(\delta)$ corresponds to assumption~(M-2) with
$(\btheta,\sigma)$
in $\mathcal{C} (\beta,R,\delta,\rho,1 )$ such that
only the first component of $\btheta$ is nonzero and $\sigma$ is
constant and
equal to one.

We write henceforth in this proof $\PP_{\theta}$ for the law
of the
process $X=(X_{t,T})_{t\leq T,T\geq1}$ and $\E_\theta$ for the corresponding
expectation.\vspace*{2pt}

Let\vspace*{2pt} $\widehat{X}=(\widehat{X}_{t,T})_{1\leq t\leq T}$ be any
predictor of
$({X}_{t,T})_{1\leq t\leq T}$ in the sense of
Definition~\ref{def:predictor}. Define
$\widehat{\theta}=(\widehat{\theta}_{t,T})_{0\leq t\leq T-1}\in\R
^{T}$ by
\[
\widehat{\theta}_{t,T}= \cases{ \widehat{X}_{t+1,T}/X_{t,T},&
\quad if $X_{t,T}\neq0$,
\cr
0, &\quad otherwise.} %
\]
For any vectors $\mathbf{u},\mathbf{v}\in\R^T$, we define
%
\begin{equation}
\label{eq:distance} d_{X}(\mathbf{u},\mathbf{v})= \Biggl(
{1\over T}\sum_{t=0}^{T-1}X_{t,T}^2
(u_t-v_t )^2 \Biggr)^{1/2}.
\end{equation}
By~(\ref{eq:lb-tvar-eq}), since $X_{t,T}$ and $\widehat{\theta
}_{t,T}$ are
$\mathcal{F}_{t,T}$-measurable, they are independent of $\xi_{t+1}$
and we have
\[
{1\over
T}\sum_{t=1}^T
\E_{\theta} \bigl[ (\widehat{X}_{t,T}-X_{t,T}
)^{2} \bigr]-1=\E_{\theta} \bigl[d_{X}^2
\bigl(\widehat\theta,\mathrm{v}_{T}\{\theta\}\bigr)
\bigr], %
\]
where, for any $\theta:(-\infty,1]\to\R$, $\mathrm{v}_{T}\{\theta\}\in\R
^{T}$ denotes the $T$-sample of $\theta$ on the regular grid
$0,1/T,\dots,(T-1)/T$,
\[
\mathrm{v}_{T}\{\theta\}= \biggl(\theta\biggl(
\frac{t}{T} \biggr) \biggr)_{0\leq
t\leq T-1} . %
\]
Hence, to prove the lower bound of
Theorem~\ref{theorem:lower_bound_quadratic_risk}, it is sufficient to
show that
there exist $\theta_{0},\ldots,\theta_{M}\in\Lambda_{1}(\beta,R)\cap
s_{1}(\delta)$, $c>0$ and $T_{0}\geq1$ both\vspace*{1pt} depending only on $\delta
$, $\beta$, $R$ and
the density $f$, such that for any
$\widehat{\theta}=(\widehat{\theta}_{t,T})_{0\leq t\leq T-1}$
adapted to
$(\mathcal{F}_{t,T})_{0\leq t\leq T-1}$ and $T\geq T_{0}$, we have
%
\begin{equation}
\label{eq:objectif:minimax} \max_{j=0,\ldots,M}\E_{\theta_{j}}
\bigl[d_{X}^2 \bigl(\widehat\theta,\mathrm{v}_{T}\{\theta_j\} \bigr) \bigr]\geq c  T^{-2\beta/ (2\beta+1 )}.
\end{equation}
We now face the more standard problem of providing a lower bound for the
minimax rate of an estimation error, since $\widehat{\theta}$ is an estimator
of $\mathrm{v}_{T}\{\theta\}$. The path for
deriving such a lower bound is
explained in
[\citeauthor{Tsybakov2009} (\citeyear{Tsybakov2009}), Chapter~2].
However, we have to
deal with a loss function $d_{X}$ which depends on the observed process~$X$.
Not only the loss function is random, but it is also not independent of the
estimator~$\widehat\theta$. The proof of the lower
bound~(\ref{eq:objectif:minimax}) thus requires nontrivial
adaptations. It
relies on some intermediate lemmas.

%
\begin{lemma}\label{lemma:birge}
We write $\mathcal{K}(\PP,\PP')$
for the Kullback--Leibler divergence between $\PP$ and $\PP'$.
For any functions $\theta_{0},\ldots,\theta_{M}$ from $[0,1]$ to $\R$
such that
%
\begin{equation}
\label{eq:kullback-condition} \max_{j=0,\ldots,M}\mathcal{K}(\PP_{\theta_{j}},
\PP_{\theta
_{0}})\leq{2\rme\over2\rme+1} \log(1+M)
\end{equation}
and any $r>0$ we have
\begin{eqnarray*}
&& \max_{j=0,\ldots,M}\E_{\theta_{j}} \bigl[d_{X}^2
\bigl(\widehat\theta,\mathrm{v}_{T}\{\theta_1\}
\bigr) \bigr]
\\
&&\qquad \geq{r^2\over4} \biggl({1\over2\rme+1}-\max
_{j=0,\ldots,M}\PP_{\theta
_{j}} \Bigl(\min_{i:i\neq j}d_{X,T}(
\theta_{i},\theta_{j})\leq r \Bigr) \biggr), %
\end{eqnarray*}
where we denote, for any two functions $\theta,\theta'$ from
$(-\infty,1]$ to $\R$,
\[
d_{X,T}\bigl(\theta,\theta'\bigr)=d_X
\bigl(\mathrm{v}_{T}\{\theta\},\mathrm{v}_{T}\{\theta'\} \bigr)  . %
\]
\end{lemma}

The proof is postponed to Section~\ref{section:proof_lemma_birge} in
\hyperref[section:useful_bounds_TVAR_postponed_proofs]{Appendix}.

We next construct certain functions $\theta_{0},\ldots,\theta_{M}\in
\Lambda_{1}(\beta,R)\cap s_{1}(\delta)$ fulfilling~(\ref
{eq:kullback-condition}) and well spread in terms of the
pseudo-distance $d_{X,T}$.
Consider the infinitely differentiable kernel $K$ defined by
\[
K(u)=\exp\biggl(-{1\over1-4u^2} \biggr)\1_{\llvert u\rrvert < 1/2} .
\]
Given\vspace*{1pt} any $m\geq8$, Vershamov--Gilbert's lemma [\citeauthor{Tsybakov2009}
(\citeyear{Tsybakov2009}), Lemma 2.9]
ensures the existence of $M+1$ points ${w^{(0)},\ldots,w^{(M)}}$ in the
hypercube $ \{0,1 \}^m$ such that
%
\begin{eqnarray}\label{eq:Vershamov}
&&M\geq2^{m/8},\qquad w^{(0)}=0\quad\mbox{and}\quad
\operatorname{card} \bigl\{ \ell: w_{\ell}^{(j)}\neq
w_{\ell}^{(i)} \bigr\}\geq m/8
\nonumber\\[-8pt]\\[-8pt]
\eqntext{\mbox{for all } j\neq i.}
\end{eqnarray}
We then define $\theta_{0},\ldots,\theta_{M}$ by setting, for all
$x\leq1$,
%
\begin{equation}
\label{eq:def-thetaj} \theta_{j}(x)={R_0\over m^{\beta}} \sum
_{\ell
=1}^{m}w_{l}^{(j)}K
\biggl(mx-\ell+\frac{1}{2} \biggr)\qquad\mbox{for } j=0,\ldots,M   ,
\end{equation}
where
%
\begin{equation}
\label{eq:L0-def} R_0 = \min\biggl(\delta,\frac{R}{ (2\llvert
K\rrvert _{\beta
} )} \biggr)
 .
\end{equation}
Since $K=0$ out of $(-1/2,1/2)$, we observe that
%
\begin{equation}
\label{eq:thetaj-zero} \theta_{j}(x)=0\qquad\mbox{for all } x\leq0 ,
\end{equation}
and
%
\begin{equation}
\label{eq:thetaj} \theta_{j}(x)={R_0\over m^{\beta}} w^{(j)}_{\lfloor mx\rfloor+1} K \biggl( \{mx \}-\frac{1}{2} \biggr)
\qquad\mbox{for all } x\in[0,1],
\end{equation}
where $ \{mx \}=mx-\lfloor mx\rfloor$ denotes the
fractional part of $mx$. Thus,
we have
%
\begin{equation}
\label{eq:thetastar} \theta^*:=\max_{0\leq j\leq M}\sup_{x\in[0,1]}
\bigl\llvert\theta_{j}(x)\bigr\rrvert\leq{R_0\rme^{-1}\over
m^{\beta}}\leq
\delta<1  .
\end{equation}
We first check that the definition of $R_0$ ensures that the $\theta
_{j}$'s are in the
expected set of parameters.

\begin{lemma}\label{lemma:admissible-part1}
For all $j=0,\ldots,M$, we have $\theta_{j}\in\Lambda_{1}(\beta,R)\cap
s_{1}(\delta)$.
\end{lemma}

The proof can be found in Section~\ref
{section:proof_lemma_admissible-part1} of \hyperref[section:useful_bounds_TVAR_postponed_proofs]{Appendix}.

Next, we provide a bound to check the required
condition~(\ref{eq:kullback-condition}) on the chosen~$\theta_j$'s.

\begin{lemma}\label{lemma:admissible-part2}
For all $j=1,\ldots,M$, we have
\[
\mathcal{K}(\PP_{\theta_{j}},\PP_{\theta_{0}})\leq\frac{8 \rme^{-2}
\kappa R_0^2}{(1-\delta^{2})\log2 }  \frac
{T}{m^{1+2\beta}} \log(1+M) , %
\]
where $\kappa$ is the constant appearing in~\textup{(I-3)}.
\end{lemma}

We prove it in Section~\ref{section:proof_lemma_admissible-part2} of
\hyperref[section:useful_bounds_TVAR_postponed_proofs]{Appendix}.

Finally, we need a control on the distances $d_{X,T}^2(\theta
_{i},\theta_{j})$.

\begin{lemma}\label{lemma:min_distance}
For any $\varepsilon>0$, there exists a constant $A$ depending only on
$\varepsilon$
and the density $f$ of $\xi$ such
that for all $m\geq16$, $T\geq4m$ and $j=0,\ldots,M$,
%
\begin{equation}
\label{eq:min-distance} \PP_{\theta_{j}} \biggl(\min_{i:i\neq j}d_{X,T}^2(
\theta_{i},\theta_{j}) \leq A {R_0^2\over m^{2\beta}}
\biggr)\leq\varepsilon+ \frac{2 R_0\rme^{-3}}{A (1-\delta)m^\beta} .
\end{equation}
\end{lemma}

The proof is postponed to Section~\ref
{section:proof_lemma_min_distance} of \hyperref[section:useful_bounds_TVAR_postponed_proofs]{Appendix}.

We can now complete the proof of
Theorem~\ref{theorem:lower_bound_quadratic_risk}.

\begin{pf*}{Proof of Theorem~\ref{theorem:lower_bound_quadratic_risk}}
Recall that $\theta_0,\dots,\theta_M$ in~(\ref{eq:def-thetaj}) are
some parameters only depending on $\beta$ and $\delta$ and a certain integer
$m\geq8$ and that, whatever the value of $m$, Lemma~\ref
{lemma:admissible-part1} insures that
$\theta_0,\dots,\theta_M$ belongs to $\Lambda_{1}(\beta,R)\cap
s_{1}(\delta)$.

Hence, it is now sufficient to show that~(\ref{eq:objectif:minimax})
holds for
a correct choice of $m$, relying on Lemmas~\ref{lemma:birge},~\ref
{lemma:admissible-part2}
and~\ref{lemma:min_distance}.
Let us set
%
\begin{equation}
\label{eq:m} m=\max\bigl\{ \bigl\lceil c_{0}T^{1/(2\beta+1)} \bigr
\rceil, 16 \bigr\}  ,
\end{equation}
where $c_0$ is a constant to be chosen. Then
$Tm^{-1-2\beta}\leq c_0^{-1-2\beta}$ and, by Lemma~\ref
{lemma:admissible-part2}, we can choose $c_0$ only depending on
$\beta$, $R$, $\kappa$ and $\delta$ so that condition~(\ref{eq:kullback-condition})
of Lemma~\ref{lemma:birge} is met. We thus get that, for any $r>0$,
\begin{eqnarray*}
&& \max_{j=0,\ldots,M}\E_{\theta_{j}} \bigl[d_{X}^2
\bigl(\widehat\theta,\mathrm{v}_{T}\{\theta_j\}
\bigr) \bigr]
\\
&&\qquad{} \geq{r^2\over4} \biggl({1\over2\rme+1}-\max
_{j=0,\ldots,M}\PP_{\theta
_{j}} \Bigl(\min_{i:i\neq j}d_{X,T}(
\theta_{i},\theta_{j})\leq r \Bigr) \biggr). %
\end{eqnarray*}
Applying\vspace*{1pt} Lemma~\ref{lemma:min_distance} with $\varepsilon=1/(4\rme
+2)$ and the previous bound with
$r^2=A R_0^2 m^{-2\beta}$, we get, as soon as $T\geq4m$,
\[
\max_{j=0,\ldots,M}\E_{\theta_{j}} \bigl[d_{X}^2
\bigl(\widehat\theta,\mathrm{v}_{T}\{\theta_j\}
\bigr) \bigr]\geq{r^2\over4} \biggl({1\over4\rme+2}-
\frac{2 R_0\rme^{-1}}{A
(1-\delta)m^\beta} \biggr)  . %
\]
The proof is concluded by observing that, as a consequence of~(\ref{eq:m}),
we can choose a constant $T_0$ only depending on
$\beta$, $R$, $\kappa$ and $\delta$ such that $T\geq T_0$ implies that
$T\geq4m$ and that the term between
parentheses is bounded by $1/(8\rme+4)$ from below.
\end{pf*}

\section{Numerical experiments}
\label{sec:numer-exper}

\begin{figure}[b]

\includegraphics{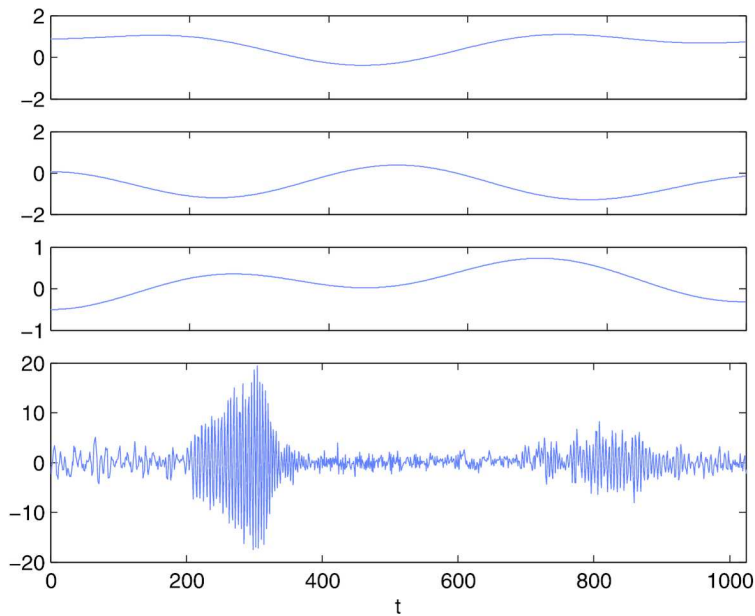}

\caption{The first three plots represent $\theta_{1}$,
$\theta_{2}$ and $\theta_{3}$ on the interval $[0,1]$. The last plot displays
$T=2^{10}$ samples of the corresponding TVAR process with Gaussian
innovations.} \label{figure:thetas_and_TVAR}
\end{figure}

In this section, we test the proposed aggregation methods on data simulated
according to a TVAR process with $d=3$. The choice of a smooth parameter
function $t\mapsto\btheta(t)$ within $s_d(\delta)$ for some $\delta
\in(0,1)$ is
done by first picking randomly some smoothly time varying partial
autocorrelation functions up to the order $d$ that are bounded between
$-1$ and
$1$ and then by relying on the Levinson--Durbin algorithm. We show the three
components of the obtained $\btheta(t)$ on $t\in[0,1]$ in the top
parts of
Figure~\ref{figure:thetas_and_TVAR}. Realizations of the TVAR process
are then
obtained from an innovation sequence $(\xi_t)_{t\in\Z}$ of i.i.d.
centered Gaussian
process with unit variance as in Definition~\ref{definition:tvar} by sampling
$\btheta$ at a given rate $T\geq1$. Figure~\ref{figure:thetas_and_TVAR}
displays one realization of such a TVAR process for $T=2^{10}$.

%

The NLMS algorithm [see \citet{GiraudRoueffSanchez-Perez2015s},
Algorithm~1] studied in
\citet{MoulinesPriouretRoueff2005} provides an online estimator of
$\btheta$
depending\vspace*{1pt} on a gradient step size $\mu$. For any $\beta\in(0,1]$, choosing
$\mu\propto T^{-2\beta/(2\beta+1)}$ yields a
$\mathcal{C} (\beta,R,\delta,\rho,1 )$-minimax-rate online
$L$-Lipschitz predictor as explained in \citet
{GiraudRoueffSanchez-Perez2015s}, Section~B.1.
Hence, proceeding as in Lemma~\ref{lemma:minimax_experts} to define
$N$ and
$\beta_i$, $i=1,\dots,N$, with $\beta_{0}=0.5$, we obtain a finite
set of NLMS predictors corresponding to gradient step sizes
$\mu_{1}>\cdots>\mu_{N}$. This set of predictors is aggregated in
two possible
ways according to the online Algorithm~\ref{algo:strat12} with the
specifications on $\eta$ and $N$ given in
Theorem~\ref{theorem:upper_bound_aggregated_forecaster_adaptative}. The
overall running time of $T$ iterates of the algorithm leading to the aggregated
predictors from the data $X_1,\dots,X_T$ is then $O(d N T)$. Since the
algorithm is recursive, the corresponding required storage capacity is
$O(d N)$.

\begin{figure}[b]

\includegraphics{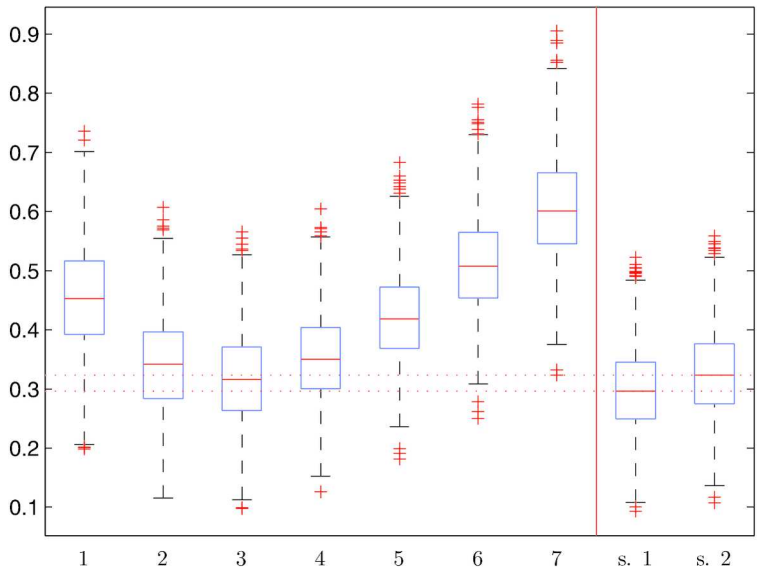}

\caption{The\vspace*{1pt} seven boxplots on the left of the vertical red line
correspond to
the averaged downward shifted empirical losses $L_{T}$ of the NLMS predictors
$\widehat{X}^{(1)},\ldots,\widehat{X}^{(7)}$. The ones on the right
of the
same line are those associated to the aggregated predictors using the
weights (\protect\ref{definition:alphas_gradient_quadratic_loss}) and
(\protect\ref{definition:alphas_quadratic_loss}).}\label
{figure:beta_0=5_T=1024}
\end{figure}

We evaluate the obtained NLMS predictors and their aggregated
predictors by
running $1000$ simulations based on equally distributed realizations of the
above Gaussian TVAR process in the case $T=2^{10}$ which yields $N=7$.
In Figure~\ref{figure:beta_0=5_T=1024}, we
compare the averaged downward shifted empirical losses defined for any
predictor $(\widehat{X}_{t,T})_{1\leq t\leq T}$ by
\[
L_T = { \frac{1}{T}}\sum
_{t=1}^{T}
\biggl( (\widehat{X}_{t,T}-X_{t,T} )^{2} -
\sigma^{2} \biggl({ \frac{t}{T}} \biggr) \biggr) .
\]
This empirical averaged loss mimics the risk considered in~(\ref{eq:max-risk}).

%

We observe that the
best NLMS predictor is the third one while the aggregated predictor of
strategy~1 enjoys a smaller loss and that of strategy~2 a slightly larger
one. This is in accordance with Theorem
\ref{thm:oracle-bounds}(i) and~(iii)
where it is
shown that the aggregated predictor of the first strategy may
outperform the
best predictor as it nearly achieves the loss of the best possible
convex combination
of the original predictors while the aggregated predictor of the second
strategy nearly achieves the loss of the best original predictor.

\begin{appendix}\label{section:useful_bounds_TVAR_postponed_proofs}
\section*{Appendix: Postponed proofs}

\subsection{A useful lemma}
The following lemma provides a uniform bound on the norm of a product of
matrices sampled from a continuous function defined on an interval $I$
and valued in a
set of $d\times d$ matrices with bounded spectral radius and norm.

%
\begin{lemma}\label{lem:product}
Let $d\geq1$ and $I$ an interval of $\R$.
Let $A$ be a function defined on $I$ taking values in the set of
$d\times
d$ matrices with eigenvalues moduli at most equal to $\delta$. Let
$\llvert \cdot\rrvert $
be any matrix norm. Denote by $A^*$ the corresponding
uniform norm of $A$,
\[
A^*=\sup_{t\in I} \bigl\vert A(t) \bigr\vert , %
\]
and, for any $h>0$, $\omega_h(A,I)$ the modulus
of continuity of $A$ over $I$,
\[
\omega_h(A;I)=\sup \bigl\{ \bigl\llvert A(t)-A(s)\bigr\rrvert :s,t
\in I, \llvert s-t\rrvert \leq h \bigr\} . %
\]
Let $\delta_1>\delta$ and assume that $A^{*}<\infty$. Then there
exist some positive constants $\varepsilon$, $\ell$ and $K$
only depending on $A^*$, $\delta$ and $\delta_1$ such that,
for any $h\in(0,1)$ fulfilling $\omega_h(A;I)\leq\varepsilon$, we
have, for all $s<t$
in $I$ and all integer $p\geq\ell(t-s)/h$,
%
\begin{equation}
\label{eq:product-univ-bound} \biggl\llvert \underbrace{A(t)A \biggl(t-\frac{t-s}{p}
\biggr)A \biggl(t-\frac
{2(t-s)}{p} \biggr)\cdots A(s)}_{p+1~\mathrm{terms}}\biggr
\rrvert \leq K \delta_1^{p+1} .
\end{equation}
\end{lemma}

\begin{pf}
Denote by $\Pi(s,t;p)$ the product of matrices appearing in the left-hand
side of~(\ref{eq:product-univ-bound}).
The proof goes along the same lines as
[Moulines, Priouret and Roueff (\citeyear
{MoulinesPriouretRoueff2005}), Proposition~13] but we use the modulus
of continuity instead of the $\beta$-Lipschitz norm to control the local
oscillation of matrices.

For $\ell_{1}\geq1$ and any square matrices $A_{1},\ldots,A_{\ell
_{1}}$, adopting the convention $\prod_{i=i_{1}}^{i_{2}}A_{i}=A_{i_{1}}\cdots A_{i_{2}}$ if $i_{1}\leq i_{2}$
and $\prod_{i=i_{1}}^{i_{2}}A_{i}$ is the identity matrix if
$i_{1}>i_{2}$, we have
%
%
\begin{eqnarray}
\label{equation:decomposition_product} \prod_{k=1}^{\ell_{1}}A_{k}
&=& A_{1}^{\ell_{1}}+\sum_{k=1}^{\ell
_{1}-1}
\Biggl(A_{1}^{\ell_{1}-k}\prod_{i=\ell_{1}-k+1}^{\ell
_{1}}A_{i}-A_{1}^{\ell_{1}-(k-1)}
\prod_{i=\ell_{1}-k+2}^{\ell
_{1}}A_{i} \Biggr)
\nonumber\\[-8pt]\\[-8pt]\nonumber
&=& A_{1}^{\ell_{1}}+\sum_{k=1}^{\ell_{1}-1}A_{1}^{\ell_{1}-k}
(A_{\ell_{1}-k+1}-A_{1} )\prod_{i=\ell_{1}-k+2}^{\ell
_{1}}A_{i}
.
\nonumber
\end{eqnarray}

Given a positive integer $\ell$, using the Euclidean division of $p+1$
by $\ell$, $p+1=\ell q+r$, we decompose the product $\Pi(s,t;p)$ as
\begin{eqnarray}\label{equation:decomposition_Pi_t_s}
\Pi(s,t;p) &=& \prod_{j=0}^{q-1} \Biggl(\prod
_{k=1}^{\ell}A \biggl(t-\frac
{(j\ell+k-1)(t-s)}{p}
\biggr) \Biggr)
\nonumber\\[-8pt]\\[-8pt]\nonumber
&&{}\times\prod_{k=1}^{r}A \biggl(t-
\frac{(q\ell+k-1)(t-s)}{p} \biggr) .
\end{eqnarray}

Using (\ref{equation:decomposition_product}), we have for any $h\geq
\ell(t-s)/p$, $0\leq j\leq q$ and $0\leq\ell_{1}\leq\ell$,
\begin{eqnarray}
\label{equation:bound_block_size_l_product}
&& \Biggl\llvert \prod_{k=1}^{\ell_{1}}A
\biggl(t-\frac{(j\ell
+k-1)(t-s)}{p} \biggr)\Biggr\rrvert
\nonumber\\[-8pt]\\[-8pt]\nonumber
&&\qquad \leq\biggl\llvert \biggl(A \biggl(t-\frac{j\ell(t-s)}{p} \biggr)
\biggr)^{\ell
_{1}}\biggr\rrvert + (\ell_{1}-1 )
\bigl(A^{*} \bigr)^{\ell
_{1}-1}\omega_{h} (A;I ) .
\end{eqnarray}

Take an arbitrary $\delta_2\in(\delta,\delta_1)$ (say the middle
point). The eigenvalues of $A$ are at most $\delta$ on $I$ and
$A^{*}<\infty$. Applying [\citeauthor{MoulinesPriouretRoueff2005}
(\citeyear{MoulinesPriouretRoueff2005}), Lemma~12] we obtain that
there is a constant
$K_{1}\geq1$ only depending on $\delta$, $\delta_{2}$ and $A^{*}$
such that $\llvert  (A(t-j\ell(t-s)/p) )^{\ell_{1}}\rrvert \leq
K_{1}\delta_{2}^{\ell_{1}}$.

From (\ref{equation:decomposition_Pi_t_s}) and (\ref
{equation:bound_block_size_l_product}), we derive the following inequality:
\[
\bigl\llvert \Pi(s,t;p)\bigr\rrvert \leq \bigl(K_1
\delta_2^\ell+ K_{2}\omega_h(A;I)
\bigr)^{q} \bigl(K_1\delta_2^r+
K_{2}\omega_h(A;I) \bigr) , %
\]
where $K_{2}= (\ell-1 )(\max\{A^{*},1\})^{\ell-1}$.

We can choose a positive integer $\ell$ and a positive number
$\varepsilon_0$
only depending on $\delta_2$, $\delta_1$ and $K_1$ such that
\[
K_1\delta_2^\ell\leq\delta_1^\ell-
\varepsilon_0 . %
\]
In the following, we set $\varepsilon=\varepsilon_0/K_{2}$. The
previous bound gives
that for any $h\in(0,1)$ such that $\omega_h(A;I)\leq\varepsilon$ and
$\ell(t-s)/p\leq h$,
\begin{eqnarray*}
\bigl\llvert \Pi(s,t;p)\bigr\rrvert &\leq& \delta_1^{\ell q}
\bigl(K_1\delta_2^r+\varepsilon_0
\bigr) \leq K_1\delta_1^{p+1}+
\varepsilon_0\delta_1^{\ell q}
\\
&\leq&
\bigl(K_1+\varepsilon_0\max \bigl\{1,
\delta_1^{1-\ell} \bigr\} \bigr) \delta_1^{p+1}. %
\end{eqnarray*}
Hence, we have the result.
\end{pf}

\subsection{Proof of Proposition~\protect\ref{proposition:stationary_TVAR}}
\label{sec:proof-prop-tvar}
We can now provide a proof of Proposition~\ref{proposition:stationary_TVAR}.

Equation~(\ref{equation:definition_TVAR}) can be more compactly
written as
%
\begin{eqnarray}
\label{equation:definition_TVARb} X_{t,T}&=&\btheta' \biggl({
\frac{t-1}{T}} \biggr)\mathbf {X}_{t-1,T}+\sigma \biggl({
\frac{t}{T}} \biggr)\xi _{t,T} .
\end{eqnarray}
%

For all $k\geq0$, iterating this recursive equation $k$ times, we have
%
\begin{eqnarray}\label{equation:TVAR_initial_condition_noise}
X_{t,T} &=& \mathbf e'_{1} \Biggl[\prod
_{i=1}^{k+1}A \biggl({ \frac{t-i}{T}} \biggr)
\Biggr]\mathbf X_{t-k-1,T}
\nonumber\\[-8pt]\\[-8pt]\nonumber
&&{} +\sum_{j=0}^{k}
\sigma \biggl({ \frac{t-j}{T}} \biggr)\mathbf e'_{1}
\Biggl[\prod_{i=1}^{j}A \biggl({
\frac
{t-i}{T}} \biggr) \Biggr]\mathbf e_{1}\xi_{t-j} ,
\end{eqnarray}
where $\mathbf e_{1}=[1\ \  0\cdots0]'$ and
\[
A(u)= \left[ \matrix{ \theta_1(u) &\theta_2(u)&\cdots& \cdots&\theta_d(u)
\cr
1 &0& \cdots& \cdots&0
\cr
0 &1&0 & \ddots&0
\cr
\vdots&0&\ddots&\ddots&\vdots
\cr
0& \cdots&0&1&0} \right].
\]

Note that the eigenvalues of $A(u)$ are the reciprocals of
the roots of the local time varying autoregressive polynomial
$z\mapsto\btheta(z;u)$, and thus are at most
$\delta<1$. Moreover, since $\btheta$ is bounded by a constant only
depending on
$d$ and is
uniformly continuous on $I=(-\infty,1]$, so is $A$ as a function
defined on $I$ and we can find
$h\in(0,1)$ such that $\omega_h(A,I)\leq\varepsilon$ for any positive
$\varepsilon$. If
$\btheta\in\Lambda_{d} (\beta,R )$, this $h$ can be
chosen depending only on $\varepsilon, \beta$ and $R$ (and also on
the matrix norm $\llvert \cdot\rrvert $).

Consider $\delta_{1}\in(\delta,1)$. Lemma~\ref{lem:product} gives
that there exist some positive constant $\varepsilon$, $\ell$ and $K$
only depending on $A^*$, $\delta$ and $\delta_1$ such that,
for any $h\in(0,1)$ fulfilling $\omega_h(A;I)\leq\varepsilon$, we
have, for all $T\geq1$, $t\leq T$ and $j\geq1$ so that $T\geq\ell/h$,
\[
\Biggl\llvert \prod_{i=1}^{j}A \biggl({
\frac
{t-i}{T}} \biggr)\Biggr\rrvert \leq K \delta_{1}^j
. %
\]
We here consider the $\ell^{\infty}$ operator norm which is the maximum
absolute row sum of the matrix, in which case $A^{*}=\max\{1,\sup_{u\in
I}(\llvert \theta_1(u)\rrvert +\cdots+\llvert \theta
_d(u)\rrvert )\}\leq2^{d}d^{1/2}$. Hence, by
(\ref{equation:TVAR_initial_condition_noise}) we obtain that
%
\begin{eqnarray}\label{eq:series-representation}
&& X_{t,T} = \sum_{i=1}^{d}
b_{t,T}(k,i) X_{t-k-i,T}+ \sum_{j=0}^{k}a_{t,T}(j)
\sigma \biggl({ \frac
{t-j}{T}} \biggr) \xi_{t-j,T},
\nonumber\\[-8pt]\\[-8pt]
\eqntext{1\leq t\leq T,}
\end{eqnarray}
with, provided that $T> \ell/h$, for all $t\leq T$, $k,j\geq1$ and
$i=1,\dots,d$,
\begin{eqnarray*}
\bigl\llvert b_{t,T}(k,i)\bigr\rrvert &\leq& K\delta_{1}^{k+1}
,
\\
\bigl\llvert a_{t,T}(j)\bigr\rrvert &\leq& K\delta_{1}^j
.
\end{eqnarray*}

The result follows.

\subsection{Proof of Lemma~\protect\ref
{lemma:bound_integral_concave_exponential}}
\label{sec:proof-lemma-technical-exp-concave}
Denote $\omega(x)=\min\{2^{-1/2},\max\{x,-2^{-1/2}\}\}$, so that
$\omega(x)^2=\min(1/2,x^2)\leq x^2$.
The\vspace*{1pt} function $x\mapsto\exp(-x^{2})$ is concave on
$[-2^{-1/2},2^{-1/2}]$, so
introducing $\omega(x)$ and then using Jensen's inequality, we get
\begin{eqnarray*}
&& \int\exp \bigl(-x^2 \bigr)\,\rmd\PP (x )
\\
&&\qquad \leq \int\exp \bigl(-
\omega^{2} (x ) \bigr)\,\rmd\PP (x ) \leq \exp \biggl(- \biggl(\int
\omega (x )\,\rmd \PP (x ) \biggr)^{2} \biggr)
\\
&&\qquad = \exp \biggl(- \biggl(\int x\,\rmd\PP (x ) \biggr)^{2}+ \biggl(\int x
\,\rmd\PP (x ) \biggr)^{2}- \biggl(\int\omega (x )\,\rmd\PP (x )
\biggr)^{2} \biggr) .
\end{eqnarray*}
It only remains to show that $(\int x\,\rmd\PP(x))^{2}-(\int
\omega(x)\,\rmd\PP(x))^{2} \leq(a^{2}-1/2)_{+}$, with the assumption that
$\PP$ has support on $[-a,a]$. This is verified if $a\leq2^{-1/2}$,
so we now assume
$a>2^{-1/2}$. We write
\begin{eqnarray*}
&& \biggl(\int x\,\rmd\PP (x ) \biggr)^{2}- \biggl(\int\omega (x )\,\rmd
\PP (x ) \biggr)^{2}
\\
&&\qquad = \int \bigl(x-\omega (x ) \bigr) \bigl(y+\omega (y
) \bigr) \,\rmd\PP (x )\,\rmd\PP (y ) .
\end{eqnarray*}
We note that $\llvert  x-\omega(x)\rrvert =(\llvert  x\rrvert -1/2)_+$ and
$\llvert  y+\omega(y)\rrvert \in\{2\llvert  y\rrvert ,\llvert  y\rrvert +2^{-1/2}\}$. We deduce that the product
$(x-\omega(x))(y+\omega(y))$ either take nonpositive values or
positive values
of the form
\[
\cases{ 2\llvert y\rrvert \bigl(\llvert x\rrvert -2^{-1/2} \bigr), &\quad
with $\llvert x\rrvert >2^{-1/2}, \llvert y\rrvert <2^{-1/2}$,
\vspace*{3pt}\cr
\bigl(\llvert x\rrvert -2^{-1/2} \bigr) \bigl(\llvert y\rrvert
+2^{-1/2} \bigr), &\quad with $\llvert x\rrvert >2^{-1/2}, \llvert
y\rrvert >2^{-1/2}$.} %
\]
Now,\vspace*{1pt} for $x,y\in[-a,a]$ with $a>2^{-1/2}$,
in the first case, we have $2\llvert  y\rrvert (\llvert  x\rrvert -2^{-1/2}) \leq2^{1/2}(a-2^{-1/2})
\leq
a^{2}-1/2$ since $2^{1/2}\leq a+2^{-1/2}$, and, in the second case,
$(\llvert  x\rrvert -2^{-1/2})(\llvert  y\rrvert +2^{-1/2})
\leq(a-2^{-1/2})(a+2^{-1/2}) =
a^{2}-1/2$. The lemma follows.

\subsection{Proof of Lemma~\protect\ref{lemma:birge}} \label
{section:proof_lemma_birge}
We define $\hat{\j}$ as the
(random) smallest index which minimizes $d_{X}(\widehat\theta,\mathrm{v}_{T}\{\theta_{j}\})$ over
$j\in\{0,\dots,M\}$ so that $d_{X}(\widehat\theta,\mathrm{v}_{T}\{\theta_{\hat{\j}}\})=\min_{\theta
\in \{\theta_{0},\ldots,\theta
_{M} \}}d_{X}(\widehat\theta,\mathrm{v}_{T}\{\theta\})$.
Note that
$d_{X,T}(\theta_{\hat{\j}},\theta_{j})\leq
d_{X}(\mathrm{v}_T\{\theta_{\hat{\j}}\},\widehat\theta)+\break d_{X}(\widehat
\theta, \mathrm{v}_{T}\{\theta_{j}\})\leq
2d_X(\widehat\theta,\mathrm{v}_T\{\theta_{j}\})$. Hence,
\begin{eqnarray*}
&& \max_{j=0,\ldots,M}\E_{\theta_{j}} \bigl[d_{X}^2
\bigl(\widehat\theta ,\mathrm{v}_{T}\{\theta_{j}\} \bigr)
\bigr]
\\
&&\qquad \geq {1\over4} \max_{j=0,\ldots,M}\E
_{\theta_{j}} \bigl[d_{X,T}^2(\theta_{\hat{\j}},
\theta_{j}) \bigr]
\\
&&\qquad \geq{r^2\over4}\max_{j=0,\ldots,M}\PP_{\theta_{j}}
\Bigl( \{ \hat{\j}\neq j \} \cap \Bigl\{\min_{i:i\neq j}d_{X,T}(
\theta _{i},\theta_{j})> r \Bigr\} \Bigr)
\\
&&\qquad \geq{r^2\over4} \Bigl(1-\min_{j=0,\ldots,M}
\PP_{\theta_{j}} (\hat{\j }= j )-\max_{j=0,\ldots,M}\PP_{\theta_{j}}
\Bigl(\min_{i:i\neq
j}d_{X,T}(\theta_{i},
\theta_{j})\leq r \Bigr) \Bigr).
\end{eqnarray*}
Birg\'e's lemma [\citeauthor{Massart2007} (\citeyear{Massart2007}),
Corollary~2.18] implies that
\[
\min _{j=0,\ldots,M}\PP_{\theta_{j}} (\hat{\j}=
j )\leq\max \biggl\{ \biggl(\frac{2\rme}{2\rme+1}\biggr),  \biggl(
\frac{\max_{j=0,\ldots,M}\mathcal{K}(\PP_{\theta_{j}},\PP _{\theta_{0}})}{\log(1+M) }\biggr) \biggr\},
\]
so the lemma follows from condition~(\ref{eq:kullback-condition}).

\subsection{Proof of Lemma~\protect\ref{lemma:admissible-part1}}
\label{section:proof_lemma_admissible-part1}
By~(\ref{eq:thetastar}), we have $\theta_j\in s_1(\delta)$ for all
$j=0,\dots,M$.
Decompose the H\"{o}lder-exponent $\beta=k+\alpha$ where $k$ is an
integer and
$\alpha\in(0,1]$. Differentiating~(\ref{eq:def-thetaj}) $k$ times,
we have, as in~(\ref{eq:thetaj}),
\[
\theta_{j}^{(k)}(x)={R_0\over m^{\alpha}} w^{(j)}_{\lfloor mx\rfloor
+1} K^{(k)} \biggl( \{mx \}-\frac{1}{2} \biggr)\qquad\mbox
{for all } x\in[0,1].
\]
Thus, for $s,s'$ in the same interval
$[\ell/m,(\ell+1)/m]$ with $\ell=0,\ldots,m-1$, we get
\begin{eqnarray*}
\bigl\llvert \theta_{j}^{(k)}(s)-\theta_{j}^{(k)}
\bigl(s'\bigr)\bigr\rrvert &\leq& {R_0\over
m^{\alpha}} \biggl
\llvert K^{(k)} \biggl(m s-\ell-\frac{1}{2} \biggr)-K^{(k)}
\biggl(m s'-\ell-\frac{1}{2} \biggr)\biggr\rrvert
\\
& \leq& R_0 \llvert K\rrvert _{\beta} \bigl\llvert
s-s'\bigr\rrvert ^\alpha.
\end{eqnarray*}
The same inequality then follows with $R_{0}$ replaced by $2R_{0}$ for
$s,s'$ in two such consecutive intervals. Now, if $s,s'$ are separated
by at least one such interval, we have
$ \vert s-s' \vert\geq m^{-1}$ and, using that $K$ has
support in $(-1/2,1/2)$, we
have that
$ \vert K^{(k)}(x) \vert$ is bounded by $ \vert K
\vert_{\beta}$. We thus get
in this case that
\[
\vert\theta_{j}^{(k)}(s)-\theta_{j}^{(k)}
\bigl(s' \bigr) \vert\leq {2R_0\over
m^{\alpha}} \sup
_{-1/2\leq x\leq1/2}\bigl\llvert K^{(k)}(x)\bigr\rrvert \leq
2R_0 \llvert K\rrvert _{\beta} \bigl\llvert s-s'
\bigr\rrvert ^\alpha . %
\]
The last two displays and~(\ref{eq:L0-def}) then yields $\theta
_{j}\in\Lambda_{1}(\beta,R)$.

\subsection{Proof of Lemma~\protect\ref{lemma:admissible-part2}}
\label{section:proof_lemma_admissible-part2}
Let $j=1,\dots,M$. Recall that $\theta_0\equiv0$ by~(\ref
{eq:Vershamov}) and~(\ref{eq:def-thetaj}).
By~(\ref{eq:thetaj-zero})
and~(\ref{eq:lb-tvar-eq}), we have that $(X_{s,T})_{s\leq0}$ has the
same distribution under $\PP_{\theta_{j}}$
and $\PP_{\theta_{0}}$ [which\vspace*{1pt} is the distribution of $(\xi_s)_{s\leq0}$].
Hence, the likelihood ratio
$\rmd\PP_{\theta_{j}}/\rmd\PP_{\theta_{0}}$ of
$(X_{s,T})_{s\leq T}$ is given by the corresponding conditional likelihood
ratio of $(X_{s,T})_{1\leq s\leq T}$ given $(X_{s,T})_{s\leq0}$. Hence,
under~(I-3), we obtain that
\begin{eqnarray*}
\frac{\rmd\PP_{\theta_{j}}}{\rmd\PP_{\theta_{0}}}&=&\prod_{t=1}^T
\frac
{f (X_{t,T}-\theta_{j}((t-1)/T)X_{t-1,T} )}{
f (X_{t,T}-\theta_{0}((t-1)/T)X_{t-1,T} )}
\\
&=& \prod_{t=1}^T \frac
{f (X_{t,T}-\theta_{j}((t-1)/T)X_{t-1,T} )}{
f (X_{t,T} )} , %
\end{eqnarray*}
where, in the second equality, we used again that $\theta_0\equiv0$.
Now, under $\PP_{\theta_j}$, we have
$X_{t,T}=\theta_{j}((t-1)/T)X_{t-1,T}+\xi_{t}$. Thus, we get
\begin{eqnarray*}
\mathcal{K} (\PP_{\theta_{j}},\PP_{\theta_{0}} ) & =& \E_{\theta_{j}}
\biggl[\log\frac{\rmd\PP_{\theta_{j}}}{\rmd\PP
_{\theta_{0}}} \biggr]
\\
&=& \sum_{t=1}^T \E_{\theta_{j}}
\biggl[\log\frac{f(\xi
_{t})}{f(\theta_{j}((t-1)/T)X_{t-1,T}+\xi_{t})} \biggr]
\\
&=& \sum_{t=1}^T \E_{\theta_{j}}\int
\log \biggl(\frac{f(u)}{f(\theta_{j}((t-1)/T)X_{t-1,T}+u)} \biggr) f(u) \,\rmd u.
\end{eqnarray*}
Using assumption~(I-3) yields
%
\begin{eqnarray}
\label{eq:etap1:admissible}
\qquad {\mathcal{K} (\PP_{\theta_{j}},\PP_{\theta_{0}} )} &\leq&\sum
_{t=1}^{T}\E_{\theta_{j}} \biggl[\kappa
\theta _{j}^{2} \biggl({ \frac{t-1}{T}}
\biggr)X_{t-1,T}^{2} \biggr] \leq\kappa\theta^{*2}\sum
_{t=1}^{T}\E_{\theta_{j}}
\bigl[X_{t-1,T}^{2} \bigr].
\end{eqnarray}
The series
representation~(\ref{eq:inf-series-representation-d1}), the fact that
$\xi$ is
centered with unit variance and~(\ref{eq:thetastar}) imply that for
all $t=0,\dots,T$
\[
\E_{\theta_j} \bigl[X_{t,T}^{2} \bigr] \leq \bigl(1-
\theta ^{*2} \bigr)^{-1} . %
\]
Using this bound and~(\ref{eq:thetastar}) in~(\ref
{eq:etap1:admissible}), we obtain
\[
{\mathcal{K} (\PP_{\theta_{j}},\PP_{\theta_{0}} )} \leq \frac{R_0^2 \rme^{-2} \kappa T}{(1-\delta^{2}) m^{2\beta}} .
\]
The proof of Lemma~\ref{lemma:admissible-part2} now follows by
applying the first bound in~(\ref{eq:Vershamov}).

\subsection{Proof of Lemma~\protect\ref{lemma:min_distance}} \label
{section:proof_lemma_min_distance}
The proof relies on an upper bound of $d_{X,T}^2(\theta_{i},\theta_{j})$
involving the noise $(\xi_{t})$. By
the expression of $\theta_j$ in~(\ref{eq:thetaj}), we have
%
\begin{equation}
\label{eq:d2XTdef-thetaj} d_{X,T}^2(\theta_{i},
\theta_{j})={ \frac
{R_0^2}{Tm^{2\beta}}}\sum_{t=0}^{T-1}X_{t,T}^{2}
\bigl(w_{k(t)}^{(i)}-w_{k(t)}^{(j)}
\bigr)^{2} K^{2} \bigl(\varphi (t ) \bigr) ,
\end{equation}
where we denoted $\varphi(t)= \{mt/T \}-1/2$ and
$k(t)=\lfloor
mt/T\rfloor+1$.
Using~(\ref{eq:inf-series-representation-d1}) and~(\ref
{eq:thetastar}), we have, for all $0\leq t\leq T-1$,
\[
\llvert X_{t,T}\rrvert \geq \llvert \xi_{t}\rrvert -\sum
_{j=1}^{\infty}\theta^{*j}\llvert \xi
_{t-j}\rrvert , %
\]
which implies
\[
X_{t,T}^2\geq \xi_{t}^2-2\llvert
\xi_{t}\rrvert \sum_{j=1}^{\infty}
\theta ^{*j}\llvert \xi_{t-j}\rrvert . %
\]
Inserting this bound in~(\ref{eq:d2XTdef-thetaj}), we get
%
\begin{equation}
\label{eq:bound-d2-remainder} \frac{m^{2\beta}}{R_0^2} d_{X,T}^2(\theta
_{i},\theta_{j})\geq\frac{1}T\sum
_{t=0}^{T-1}\xi_{t}^2
\bigl(w_{k(t)}^{(i)}-w_{k(t)}^{(j)}
\bigr)^{2} K^{2} \bigl(\varphi (t ) \bigr) -
\mathcal{R}_T ,
\end{equation}
where
\[
\mathcal{R}_T=\frac{2\rme^{-2}}{T}\sum_{t=0}^{T-1}
\sum_{j=1}^{\infty}\theta^{*j}\llvert
\xi_{t}\rrvert \llvert \xi _{t-j}\rrvert. %
\]
Thus, with~(\ref{eq:bound-d2-remainder}), the left-hand side of
inequality~(\ref{eq:min-distance}) is upper bounded by
\[
\PP_{\theta_j} \Biggl(\min_{i:i\neq j}\frac{1}T\sum
_{t=0}^{T-1}\xi _{t}^2
\bigl(w_{k(t)}^{(i)}-w_{k(t)}^{(j)}
\bigr)^{2} K^{2} \bigl(\varphi (t ) \bigr)< 2A \Biggr) + \PP(
\mathcal{R}_T>A) .
\]
Using that $\xi$ is centered with unit variance and then~(\ref
{eq:thetastar}), we easily get that
\[
\E_{\theta_j} [\mathcal{R}_T ]\leq\frac{2\rme
^{-2}}T
\sum_{t=0}^{T-1}\sum
_{j=1}^{\infty}\theta^{*j} \leq
\frac{2\rme^{-2}\theta^*}{1-\theta^*} \leq\frac{2 R_0\rme
^{-3}}{(1-\delta)m^\beta} .
\]
Hence, by Markov's inequality, to conclude the proof, it now suffices
to show
that, for $A$ well chosen,
%
\begin{equation}
\label{eq:bound-d2-remainderRT} \PP_{\theta_j} \Biggl(\min_{i:i\neq j}
\frac{1}T\sum_{t=0}^{T-1}\xi
_{t}^2 \bigl(w_{k(t)}^{(i)}-w_{k(t)}^{(j)}
\bigr)^{2} K^{2} \bigl(\varphi (t ) \bigr)< 2A \Biggr) \leq
\varepsilon .
\end{equation}
For $k\in \{1,\ldots,m \}$ we define
$J_{k}= \{\lfloor(k-1)T/m\rfloor+i  \lceil T/(4m)\rceil+1
\leq i \leq\lfloor3T/(4m)\rfloor \}$.
We observe that the cardinality of $J_{k}$ is
\[
\Gamma \biggl(\frac{T}{m} \biggr)= \biggl\lfloor\frac{3T}{4m} \biggr
\rfloor- \biggl\lceil\frac{T}{4m} \biggr\rceil\geq1 ,
\]
where the lower bound is a consequence of the assumption $T\geq4 m$ in
the lemma.
Moreover, it is easy to check that we have $\vert\varphi(t)
\vert\leq1/4$ for
all index
$t\in J_{k}$ and that, for each $1\leq k\leq m$, the set $J_k$ is
included in the set $\{1\leq t\leq T-1:k(t)=k\}$ (so
that, in particular,
$J_k\cap J_{k'}=\varnothing$ for $k<k'$). It follows that random variables
\[
S_{k}={1\over\Gamma(T/m)} \sum_{t\in J_{k}}
\xi_{t-1}^{2}\qquad \mbox{for } k=1,\ldots,m
\]
are i.i.d. By the monotonicity of $K$ in $\R_{-}$ and its symmetry, we have
\begin{eqnarray*}
\frac{1}T\sum_{t=0}^{T-1}
\xi_{t}^2 \bigl(w_{k(t)}^{(i)}-w_{k(t)}^{(j)}
\bigr)^{2} K^{2} \bigl(\varphi (t ) \bigr) &\geq&
{1\over T}\sum_{k=1}^{m}
\bigl(w_{k}^{(i)}-w_{k}^{(j)}
\bigr)^{2}\sum_{t\in J_k}\xi_{t}^{2}K^{2}
\bigl(\varphi (t ) \bigr)
\\
&\geq& {K^2(1/4)\Gamma(T/m) \over T}\sum_{k=1}^{m}
\bigl(w_{k}^{(i)}-w_{k}^{(j)}
\bigr)^{2}S_{k}.
\end{eqnarray*}
From~(\ref{eq:Vershamov}), for any $i,j\in \{1,\ldots,M \}
$ there exist
at
least $\lceil m/8\rceil$ values of $k$ for which
$(w_{k}^{(i)}-w_{k}^{(j)})^{2}$ equals one in the
above sum. Hence, using the order statistics $S_{(1,m)}\leq\cdots\leq
S_{(m,m)}$, we thus obtain that
\begin{eqnarray*}
\min_{i:i\neq j}{1\over T}\sum
_{t=0}^{T-1}\xi_{t}^{2}
\bigl(w_{k(t)}^{(i)}-w_{k(t)}^{(j)}
\bigr)^{2}K^{2} \bigl(\varphi (t ) \bigr) &\geq&
{K^2(1/4)\Gamma(T/m)\over T}\sum_{k=1}^{\lceil m/8\rceil
}S_{(k,m)}
\\
&\geq& {{K^2(1/4) m \Gamma(T/m)}\over{16 T}} S_{(\lfloor
m/16\rfloor,m)}
\\
&\geq&{{{K^2(1/4)}\over128}} S_{(\lfloor m/16\rfloor,m)},
\end{eqnarray*}
where we used $\Gamma(T/m)\geq T/(8m)$ for $T/m\geq4$ in the last inequality.
Let us denote by $F$ the cumulative distribution function of $S_{1}$, which
only depends on $\Gamma(T/m)$ and on the distribution of $\xi_{0}$.
For $x>0$, we have
\begin{eqnarray*}
\PP (S_{(\lfloor m/16\rfloor,m)}\leq x )&=&\PP \biggl(\operatorname{Bin}\bigl(m,F(x)\bigr)
\geq \biggl\lfloor\frac{m}{16} \biggr\rfloor \biggr)
\\
&\leq& {m\over\lfloor m/16\rfloor} F(x)\leq32 F(x).
\end{eqnarray*}
Gathering the last two bounds, we get that
\begin{eqnarray*}
&& \PP_{\theta_{j}} \Biggl(\min_{i:i\neq j}{1\over T}
\sum_{t=1}^{T-1}\xi_{t}^{2}
\bigl(w_{k(t)}^{(i)}-w_{k(t)}^{(j)}
\bigr)^{2}K^{2} \bigl(\varphi (t ) \bigr)\leq2A \Biggr)
\\
&&\qquad \leq
\PP \biggl(S_{(\lfloor m/16\rfloor,m)}\leq{256 A\over
K^2(1/4)} \biggr)
\\
&&\qquad \leq 32 F \biggl(\frac{256 A}{K^2(1/4)} \biggr).
\end{eqnarray*}
Recall that $\Gamma(T/m)\geq1$ and note that $S_1$ admits a density, since
$\xi$ does. By the strong law of large numbers, we further have that the
random variable $S_{1}$ converges to
$1$ almost surely when $\Gamma(T/m)$ goes to infinity, so there exists $x_0>0$
depending only on the density of $\xi$ such that
$F(x_0)\leq\varepsilon/32$ whatever the value of \mbox{$\Gamma(T/m)\geq
1$}. Therefore,
there exists some $A>0$, depending only on the distribution of $\xi$,
such that~(\ref{eq:bound-d2-remainderRT}) holds, which achieves the proof.
\end{appendix}

\section*{Acknowledgements}
We gratefully acknowledge the fruitful comments of the referees.

\begin{supplement}[id=suppA]
\stitle{Supplementary material for: Aggregation of predictors for nonstationary sub-linear processes and online adaptive forecasting of time
varying autoregressive processes}
\slink[doi]{10.1214/15-AOS1345SUPP} 
\sdatatype{.pdf}
\sfilename{AOS1345\_supp.pdf}
\sdescription{We explain how to build nonadaptive minimax
predictors which can be used in the aggregation step. The document also
contains some technical proofs and provides additional results with
improved aggregation rates.}
\end{supplement}

%

\printaddresses

\end{document}